\documentclass[11pt]{article}
\usepackage[dvips]{graphics}
\input{amssym.def}
\input{amssym.tex}

\title{\Large\bf  Slit maps in the study of equal-strength cavities in 
$n$-connected elastic planar domains}
\author{\bf {\sc Y.A. Antipov 
}\\ 
Department of Mathematics, Louisiana State University\\
Baton Rouge LA 70803, USA}

\date{}

\newcommand{\I}{\mathop{\rm Im}\nolimits}
\newcommand{\R}{
\mathop{\rm Re}\nolimits}

\newcommand{\ov}[1]{\overline{#1}}

\newcommand{\Ga}{\alpha}
\newcommand{\Gb}{\beta}
\newcommand{\Gd}{\delta}

\newcommand{\Gve}{\varepsilon}
\newcommand{\Gf}{\phi}

\newcommand{\Gg}{\gamma}

\newcommand{\Gk}{\kappa}

\newcommand{\Gl}{\lambda}
\newcommand{\Gn}{\eta}
\newcommand{\Gm}{\mu}

\newcommand{\Gt}{\theta}

\newcommand{\Gr}{\rho}

\newcommand{\Gs}{\sigma}

\newcommand{\Go}{\omega}
\newcommand{\Gx}{\xi}

\newcommand{\Gz}{\zeta}
\newcommand{\GD}{\Delta}
\newcommand{\GF}{\Phi}
\newcommand{\GG}{\Gamma}

\newcommand{\GY}{\Psi}

\newcommand{\CD}{{\cal D}}

\newcommand{\CI}{{\cal I}}
\newcommand{\CJ}{{\cal J}}

\newcommand{\CR}{{\cal R}}

\newcommand{\beq}{\begin{equation}}
\newcommand{\eeq}{\end{equation}}
\newcommand{\barr}{\begin{eqnarray}}
\newcommand{\earr}{\end{eqnarray}}
\newcommand{\beqn}{\begin{equation*}}
\newcommand{\eeqn}{\end{equation*}}
\newcommand{\barrn}{\begin{eqnarray*}}
\newcommand{\earrn}{\end{eqnarray*}}

\newcommand{\fr}{\frac}

\begin{document}
\maketitle

\begin{abstract}

The inverse problem of plane elasticity on $n$ equal-strength cavities in a plane
subjected to constant loading at infinity and in the cavities boundary is analyzed.
By reducing the governing boundary value problem to the Riemann-Hilbert 
problem on a symmetric Riemann surface of genus $n-1$ a family of conformal mappings
from a parametric slit domain  onto the $n$-connected elastic domain is constructed.
The conformal mappings  are presented in terms of hyperelliptic integrals and the zeros
of the  first derivative of the mappings are analyzed. It is shown that for any $n\ge 1$ there always exists a set of
the loading parameters  for which these zeros generate inadmissible poles of the solution.

\end{abstract}

\setcounter{equation}{0}

\section{Introduction}

Analysis of the stresses induced by the presence of inclusions and cavities in an elastic
matrix subjected to loading have been of interest for over a century \cite{mus}, \cite{sav}, \cite{kos}. 
Inverse problems of elasticity which concern problems of determination of the shapes of curvilinear  inclusions
and cavities with prescribed properties  excite particular attention due to their relevance
to the material design  \cite{neu},   \cite{che}, \cite{vig1}, \cite{ban}, \cite {shi}.
A considerable amount of work examines equal-strain inclusions subjected to uniform loading
with uniform distribution of stresses inside. By analyzing the stress distribution 
in composites with single elliptic and ellipsoidal inclusions in two- and three-dimensional unbounded 
elastic bodies Eshelby \cite{esh} established that the stress fields are uniform in the interior 
of the inclusions provided the matrix is loaded uniformly at infinity. He also conjectured that 
there do not exist other shapes of single inclusions with such a property. This conjecture was proved for
 the plane and anti-plane model problem in the simply-connected case  in  \cite{sen}. An alternative proof for the antiplane
case by the method of conformal mappings was proposed in \cite{ru}.

Motivated by the  problem of designing perforated  structures of minimum weight Cherepanov  
\cite{che} studied  the inverse problem of elasticity on a plane uniformly
loaded at infinity and having $n$ holes. The boundary of the holes  
is subjected to constant normal and tangential traction, and the holes profile $L_j$
($j=0,1,\dots,n-1$) are determined from an extra boundary condition. It states that
the tangential normal stress $\Gs_t$ is the same constant, $\Gs$,  in all the 
contours $L_j$. 
For the solution, a conformal map of an $n$-connected  slit domain $\CD^e$
into the elastic domain $D^e$, the exterior of the $n$ holes,  is applied. The map
transforms the boundary value problem into two Schwarz problems of the theory of analytic functions 
on the $n$ slits. The feature of the map $z=\Go(\Gz)$
employed \cite{che}
is that it maps the  point $\Gz=\infty$ into the point $z=\infty$, and the exterior of $n$ parallel slits
of the parametric plane into the $n$-connected domain $D^e$.
In general,  unless $n\le 3$, for such a map
these slits do not lie in the same line. 
Cherepanov  \cite{che} solved by quadratures the problem 
in the simply and doubly connected symmetric cases and also analyzed the periodic and doubly periodic
problems. Vigdergauz \cite{vig1} noticed and corrected an error in the computations \cite{che}
implemented for the symmetric doubly connected case when the slit domain is the exterior of the cuts 
$[-2,-1]$, $[1,2]$, and the map $z=\Go(\Gz)$ has the property $\Go:\infty\to\infty$.
For the symmetric case \cite{che}, the loading parameters $a$ and $b$ are real, and
the two equal-strength holes exist if $|b/a|<1$ \cite{vig1}.

To solve the Cherepanov problem for any $n$-connected domain, Vigdergauz \cite{vig1} proposed
to employ a circular map from the exterior of $n$-circles onto the $n$-connected elastic domain.
Application of the Sherman integral representation reduced the resulting boundary value problem
to integral equations solved numerically by the method of least squares.
This method was further developed  in \cite{vig2} for doubly periodic structures by using an
integral representation with quasi-automorphic analogues of the Cauchy kernel and the numerical method
of least squares for solving the governing integral equations. 
An alternative explicit representation in terms of the Weierstrass elliptic function for the profile
of an inclusion in the case of a doubly periodic structure was given in \cite{gra}.
The antiplane shear problem for two equal-strain inclusions by means of the Weierstrass zeta function
was treated in \cite{kan}. 

The main goal of this paper was to derive a closed-form representation of a family of conformal mappings
$z=\Go(\Gz)$, $\Go: \CD^e\to D^e$ for $n= 2$  in the non symmetric case and for $n\ge 3$ and therefore to
determine a family of the profiles of $n$ equal-strength holes. In addition, we aim to derive necessary and sufficient conditions for the solution to exist.
This will be achieved by studying the
poles of one of the Kolosov-Muskhelishvili complex potentials due to possible zeros of the 
derivative $\Go'(\Gz)$ of the conformal mapping. The stress field $\{\Gs_1, \Gs_2, \tau_{12}\}$
is expressible through two functions $\GF$ and $\Psi$ \cite{mus} which have to be analytic
everywhere in the domain $D^e$. One of them, $\GF(z)$, is a constant, while the second one
has the form  \cite{che}, p. 918
\beq
\Psi(\Go(\Gz))=\fr{F_+(\Gz)+F_-(\Gz)}{2\Go'(\Gz)},\quad \Gz\in\CD^e={\Bbb C}\setminus l,
\label{1.1}
\eeq
where $F_\pm(\Gz)$, the solutions to certain Schwarz problems,  
are analytic  functions in ${\Bbb C}\setminus l$, and $l$ is the union of the slits $l_j$, $j=0,1,\ldots, n-1$, 
in the parametric $\Gz$-plane.
If $\Go'(\Gz)$ has zeros in the slit domain $\CD^e$, then the potential $\Psi(z)$ has inadmissible poles
at the images of these zeros. In \cite{vig4}, under the assumption that the solution $\Psi(z)$ found is analytic in the exterior 
of $n$ holes by applying the maximum principle it was shown that the condition 
$|(\Gs_2^\infty-\Gs_1^\infty)/(\Gs_2^\infty+\Gs_1^\infty)|\le 1$ is necessary for the existence of the solution. Here, $\Gs_1^\infty$ and $\Gs_2^\infty$ are constant stresses applied at infinity.
 To the author's knowledge, no sufficient solvability conditions for the case $n\ge 3$ and for
two nonsymmetric cavities and associated exact representations for equal-strength cavities profiles are available in the literature.

In Section 2, we formulate the problem as two Schwarz problems for two auxiliary functions coupled by two conditions. These conditions guarantee that the potential $\Psi$ is analytic everywhere in the domain $D^e$
and that the conformal map is single valued.
Section 3 gives an integral representation  in terms of elliptic integrals of the mapping function for $n=2$
in the general not necessarily symmetric case. The map has two free parameters and, in addition,
has two free scaling parameters. It is shown that if $\Gg=|b/a|<1$, then the solution
always exists. 
Here,
\beq
 \quad a=\fr12(\Gs-p)+i\tau,\quad b=\fr12(\Gs_2^\infty-\Gs_1^\infty)+i\tau^\infty,
\label{1.3}
\eeq
$\Gs_1^\infty$, $\Gs_2^\infty$, 
 and $\tau^\infty$ are constant stresses applied at infinity,  $p$
 and $\tau$ are the traction components applied to the holes boundaries, and $\Gs=\Gs_1^\infty+\Gs_2^\infty-p$.
If $\Gg>1$, then the function $\Psi(\Gz)$ has  four poles, while the contours $L_1$ and $L_2$, the profiles of the holes, may  (for sufficiently large values of the  parameter $\Gg$) or may not intersect.
If $\Gg=1$, then  the contours are two straight segments, and the function $\Psi(z)$ has removable
 singularities on the contours.
In Section 4, we analyze the triply connected case. To pursue the goal to describe the general
family of equal-strength  cavities, we  construct the most general 
form of the conformal map possible in the case $n=3$. It maps the exterior of three slits $[-1/k,-1]$, $[k_1,k_2]$, and $[1,1/k]$  ($0<k<1$, $-1<k_1<k_2<1$)
into the triply connected domain $D^e$, while  the infinite point $z=\infty$ is the image of a point $\Gz_\infty=\Gz_\infty'+i\Gz''_\infty$. Thus, in addition to  two scaling parameters it has five real free parameters, 
$k$, $k_1$, $k_2$, $\Gz_\infty'$, and $\Gz_\infty''$. The governing boundary value problem
reduces to a symmetric Riemann-Hilbert problem on a genus-2 Riemann surface
and similarly to \cite{ant} is solved exactly.
It is shown that if $\Gz_\infty$ is a finite point, then regardless of the values of the parameter $\Gg$  
the function $\Psi(z)$ has always inadmissible poles in $D^e$, and the solution does not exist, that is
the condition $\Gg<1$ is necessary but not sufficient.
In Section 5, we identify a family of $n$-connected domains which can be interpreted as images
of a slit domain with cuts locating in the same line. 
 We choose
the parametric domain $\CD^e$
as the union of slits in the real axis and assume that $\Go:\infty\to\infty$. The conformal map with such properties has $2n-2$
free real parameters.        
In particular, for the case $n=3$ and $\Gz_\infty=\infty$, we show that the family of conformal mappings
is four-parametric (the slits are $[-1, k_1]$, $[k_2,k_3]$, and $[k_4,1]$), and if $\Gg<1$, then
the solution is free of poles. Otherwise, if $\Gg>1$, it has six poles, and the solution does not exist.
In Appendix we analyze the case $n=1$ and show that the potential $\Psi(z)$ has two inadmissible poles if $\Gg>1$. The same solvability condition for the case $n=1$ by a different method was derived earlier in \cite{vig3}.

\setcounter{equation}{0}

\section{Formulation}\label{s2.1}

\begin{figure}[t]
\centerline{
\scalebox{0.5}{\includegraphics{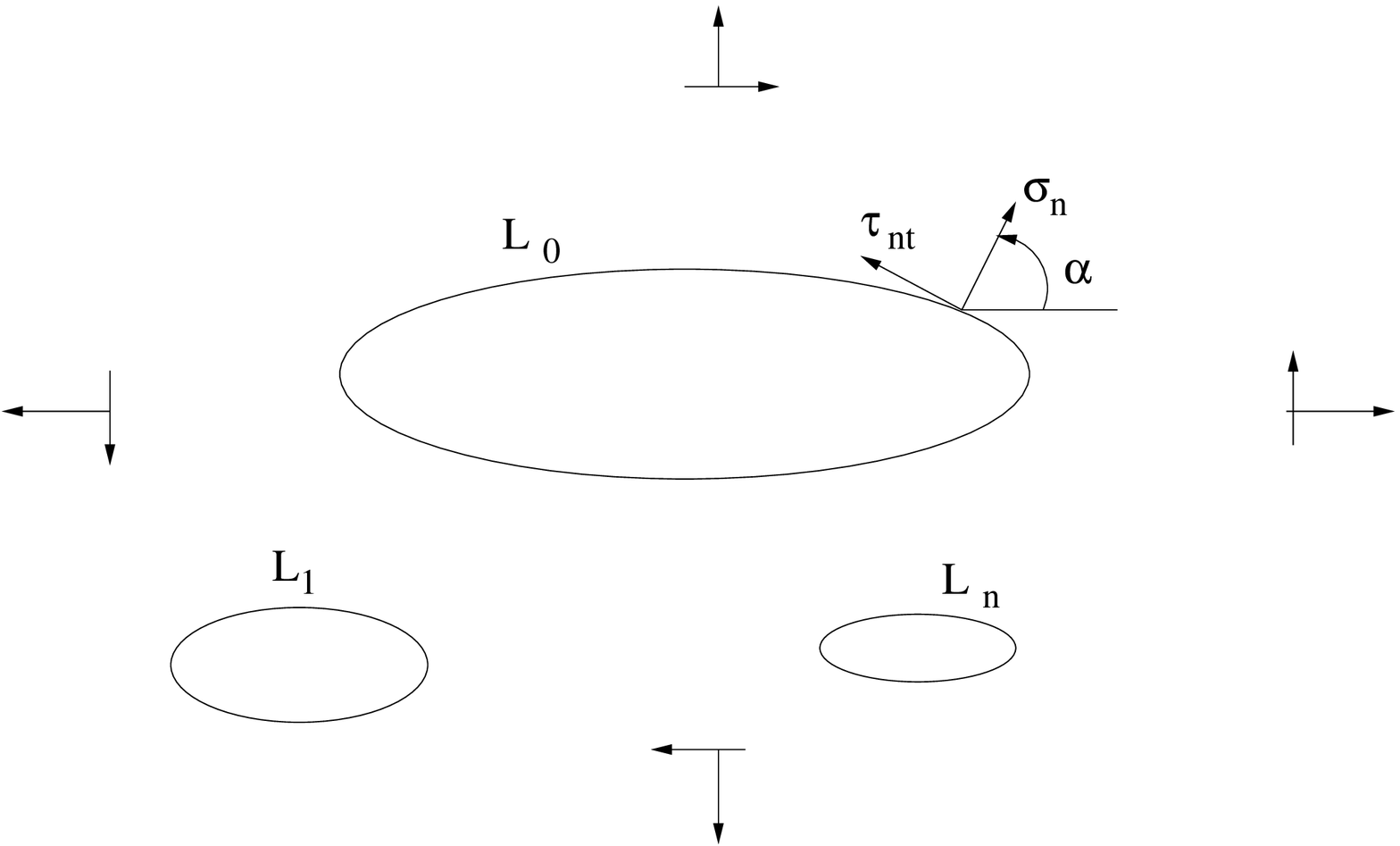}}}
\caption{
Geometry of the problem}
\label{fig1}
\end{figure}

Consider the following problem of plane elasticity  \cite{che} (Figure \ref{fig1}).

{\sl Let an infinite isotropic plane subjected to constant stresses at infinity, $\Gs_1=\Gs_1^\infty$,
$\Gs_2=\Gs_2^\infty$, and $\tau_{12}=\tau^\infty$, have $n$ holes $D_0$, $D_1$, $\ldots D_{n-1}$.
Assume that constant normal and tangential traction components are applied to their boundaries
$L_j$, $\Gs_n=p$, $\tau_{nt}=\tau$, $j=0,1,\ldots,n-1$.
Find the shape and location of the holes  such that the tangent normal stress $\Gs_t$ is constant, $\Gs_t=\Gs$,
in all the contours $L_j$.}

Let $\GF(z)$ and $\Psi(z)$ $(z=x_1+ix_2$) be the Kolosov-Muskhelishvili  potentials of the problem. These functions are analytic everywhere in the $n$-connected domain $ D^e={\Bbb C}\setminus D$, $D=\cup_{j=0}^{n-1} D_j$ and continuous
in $D^e$ up to its boundary.
The equilibrium equations of plane elasticity are satisfied if the stresses are  \cite{mus} 
\beq
\Gs_1+\Gs_2=4\R\GF(z),
\quad
\Gs_2-\Gs_1+2i\tau_{12}=2[\bar z\GF'(z)+\Psi(z)].
\label{2.1}
\eeq
The stresses  and the traction vector components on the boundaries are connected by the relations
\beq
\Gs_t+\Gs_n=\Gs_1+\Gs_2,
\quad
\Gs_t-\Gs_n+2i\tau_{nt}=e^{2i\Ga(z)}(\Gs_2-\Gs_1+2i\tau_{12}).
\label{2.2}
\eeq
Here,  $\Ga(z)$ is the angle between the positive direction of the $x_1$-axis and the external normal $n$ to
the cavity boundary (internal with respect to the body $D^e$). At infinity, the functions $\GF(z)$ and $\Psi(z)$ behave as  \cite{mus} 
$$
\GF(z)=b'+\fr{X+iY}{2\pi(1+\Gk)}\fr{1}{z}+O\left(\fr{1}{z^2}\right),
$$
\beq
\GY(z)=b-\fr{\Gk(X-iY)}{2\pi(1+\Gk)}\fr{1}{z}+O\left(\fr{1}{z^2}\right),
\label{2.2'}
\eeq
where 
$$
b'=\fr{\Gs_1^\infty+\Gs_2^\infty}{4}+iC',\quad 
b=\fr12(\Gs_2^\infty-\Gs_1^\infty)+i\tau^\infty, \quad \Gk=\fr{\Gl+3\mu}{\Gl+\mu}, 
$$
\beq
 X+iY=\sum_{j=0}^{n-1}(X_j+iY_j),\quad  X_j+iY_j=\int_{L_j}(X_{jn}+iY_{jn})ds,
\label{2.2''}
\eeq
$\Gl$ and $\mu$ are the Lam\'e constants,  $(X_j,Y_j)$  is the total force applied to the boundary $L_j$
$(j=0,\ldots,n-1)$ from the side of a normal directed towards the body $D^e$, and $X_{jn}$ and $Y_{jn}$ are the $x_1$-
and $x_2$-projections of the force at a point of the boundary.
Since the traction components 
$\Gs_n$ and $\tau_{nt}$ are constant on the boundary, $X=Y=0$.

The boundary condition $\R\GF(z)=\fr14(\Gs_1^\infty+\Gs_2^\infty)$ on $L_j$, $j=0,\ldots,n-1$, and the condition at infinity (\ref{2.2'}) 
imply  that the function $\GF(z)$ is a constant everywhere in $D^e$,
$\GF(z)= \fr14(\Gs_1^\infty+\Gs_2^\infty)+iC'$. The constant $C'$ does not affect the stresses. It is expressed
through the rotation $\Gve_\infty$ of the plane at infinity, $C'=2\mu(1+\Gk)^{-1} \Gve_\infty$  \cite{mus}.

Due to (\ref{2.1}), (\ref{2.2}), (\ref{2.2'}) the analytic function $\Psi(z)$ has to satisfy the conditions
\beq
\Psi(z)=ae^{-2i\Ga(z)}, \quad z\in L=\cup_{j=0}^{n-1}L_j, \quad 
\Psi(z)=b+O(z^{-2}), \quad z\to\infty, 
\quad a=\fr12(\Gs-p)+i\tau.
\label{2.4}
\eeq
It is known   \cite{cou}, \cite{kel} that there exists an analytic function $z=\Go(\Gz)$
that conformally maps the extended complex $\Gz$-plane ${\Bbb C}\cup\infty$ cut along $n$ segments
parallel to the real $\Gz$-axis onto the $n$-connected domain $D^e$ in the $z$-plane.
Such a function is a one-to-one map. The infinite point $z=\infty$ is the image of a certain
point $\Gz=\Gz_\infty$, and in the vicinity of that point, the conformal map $\Go(\Gz)$ can be represented as
\beq
\Go(\Gz)=\fr{c_{-1}}{\Gz-\Gz_\infty}+c_0+\sum_{j=1}^\infty c_j(\Gz-\Gz_\infty)^j.
\label{2.5}
\eeq
If $\Gz_\infty=\infty$, then 
\beq
\Go(\Gz)=c_{-1}\Gz+c_0+\sum_{j=1}^\infty \fr{c_j}{\Gz^j}.
\label{2.5'}
\eeq

The boundary condition in (\ref{2.4}) written in the $\Gz$-plane reads  \cite{che}
\beq
\psi(\Gz)\Go'(\Gz)+a\ov{\Go'(\Gz)}=0, \quad \Gz\in l=\cup_{j=0}^{n-1}l_j,
\label{2.6}
\eeq
where $\psi(\Gz)=\Psi(z)$. 
The problem is significantly simplified for the two analytic functions 
\beq
F_\pm(\Gz)=[\psi(\Gz)\pm\bar a]\Go'(\Gz).
\label{2.7.0}
\eeq
The new functions solve the following  boundary value problem.

{\sl Find two functions $F_\pm(\Gz)$ analytic in the domain ${\Bbb C}\setminus l$ continues up to the boundary $l$,
having at most integrable singularities at the endpoints of the contour $l$ and satisfying the 
Schwarz boundary conditions on the cuts $l$
\beq
\R F_+(\Gz)=0, \quad \I F_-(\Gz)=0, \quad \Gz\in l.
\label{2.7}
\eeq
If $\Gz_\infty$ is a finite point, then the functions  $F_\pm(\Gz)$ behave at $\Gz_\infty$ and $\infty$ as
\beq
F_\pm(\Gz)=-\fr{c_{-1}(b\pm \bar a)}{(\Gz-\Gz_\infty)^2}+O(1), \quad \Gz\to\Gz_\infty,
\quad 
F_\pm(\Gz)=O(\Gz^{-2}), \quad \Gz\to\infty.
\label{2.7'}
\eeq
In the case $\Gz_\infty=\infty$,
\beq
F_\pm(\Gz)=c_{-1}(b\pm \bar a)+O(\Gz^{-2}),  \quad \Gz\to\infty.
\label{2.7''}
\eeq
The two Schwarz  problems are coupled by the conditions

(i) the difference $F_+(\Gz)-F_-(\Gz)$ may not have zeros in the domain  ${\Bbb C}\setminus l$;

(ii) the following $n$ integrals over the loops $l_j$ vanish:
\beq
\int_{l_j}[F_+(\Gz)-F_-(\Gz)]d\Gz=0, \quad j=0,1,\ldots n-1.
\label{2.7'''}
\eeq
}

If the two functions $F_+(\Gz)$ and $F_-(\Gz)$ are known, then  the functions $\Go'(\Gz)$ and $\psi(\Gz)$ are determined by
\beq
\Go'(\Gz)=\fr{F_+(\Gz)-F_-(\Gz)}{2\bar a}, \quad \psi(\Gz)=\fr{F_+(\Gz)+F_-(\Gz)}{2\Go'(\Gz)}, 
\quad \Gz\in {\Bbb C}\setminus l.
\label{2.8}
\eeq
The condition (i)  is necessary and sufficient for  the function $\psi(\Gz)$ to be analytic in the domain ${\Bbb C}\setminus l$,
while the condition (ii) guarantees that the function $\Go(\Gz)$ is single-valued. 

To justify the relations (\ref{2.7'}) in the case $|\Gz_\infty|<\infty$, we note that  $\Gz=\infty$ is the $\Go$-image of a finite point of 
the domain $D^e$,
 the function $\psi(\Gz)$ is bounded
at infinity and also
\beq
\Go(\Gz)=d_0+\fr{d_1}{\Gz}+\fr{d_2}{\Gz^2}+\ldots, \quad \Go'(\Gz)=-\fr{d_1}{\Gz^2}-\fr{2d_2}{\Gz^3}-\ldots,
\quad \Gz\to\infty.
\label{2.21.0}  
 \eeq
 In a neighborhood of the point $\Gz_\infty\in(-1,1)$ the function $\Go(\Gz)$
 admits the expansion (\ref{2.5}) and therefore 
 \beq
 \psi(\Gz)=b+O((\Gz-\Gz_\infty)^2), \quad \Go'(\Gz)= -\fr{c_{-1}}{(\Gz-\Gz_\infty)^2}+ c_1+
 2 c_2(\Gz-\Gz_\infty)+\ldots, \quad \Gz\to\Gz_\infty.
 \label{2.21.1}
 \eeq

\setcounter{equation}{0}
\section{Two cavities}\label{s2.4}

\subsection{Two cavities when $\Gz_\infty$ is a finite point: the general case}
  
Every doubly connected domain $D^e$ may be conformally mapped by a function $\Gz=\Go^{-1}(z)$
 onto a slit domain $\CD^e$, the extended complex $\Gz$-plane cut along the two cuts $l_0=[-1/k,-1]$ and $l_1=[1,1/k]$,
 where $k\in(0,1)$. Moreover, it is possible to choose the map such that the infinite point $z=\infty$ falls into 
 a point $\Gz_\infty$ in the open segment $(-1,1)$ of the real $\Gz$-axis.
 Such a map can be expressed through  elliptic integrals.

We introduce  an elliptic surface $\CR$ of the algebraic function  $u^2=p_2(\Gz)$, where
$p_2(\Gz)=(\Gz^2-1)(\Gz^2-1/k^2)$. A single branch $f(\Gz)$ of the function $p_2^{1/2}(\Gz)$ is fixed 
in the $\Gz$-plane cut along the two segments $l_0$  and $l_1$ 
by the condition $p_2^{1/2}(\Gz)\sim \Gz^2$, $\Gz\to \infty$. This branch 
is pure imaginary on the sides of the cuts,
$f^\pm(\Gz)=\mp i(-1)^j\sqrt{|p_2(\Gz)|}$,  $\Gz\in l_j$, $j=0,1$, and real
for $\Gz=\Gx$ lying outside the cuts in the real axis, $f(\Gx)>0$, $|\Gx|>1/k$, and  $f(\Gx)<0$, $-1<\Gx<1$.

Since $\R F_+(\Gz)=0$ and $\I F_-(\Gz)=0$ on $l$, the functions $iF_+(\Gz)$ and $F_-(\Gz)$
 can analytically and symmetrically be continued onto the whole Riemann surface. 
 The new functions, $F_1(\Gz,u)$  and $F_2(\Gz,u)$, given by
 \beq
 F_1(\Gz,u)=\left\{ 
 \begin{array}{cc}
 iF_+(\Gz), & (\Gz,u)\in {\Bbb C}_1,\\
 -i\ov{F_+(\bar \Gz)}, & (\Gz,u)\in {\Bbb C}_2,\\
\end{array}\right. \quad 
 F_2(\Gz,u)=\left\{ 
 \begin{array}{cc}
 F_-(\Gz), & (\Gz,u)\in {\Bbb C}_1,\\
 \ov{F_-(\bar \Gz)}, & (\Gz,u)\in {\Bbb C}_2,\\
\end{array}\right. 
\label{2.0}
\eeq
 are rational on the surface and
 satisfy the symmetry condition
 \beq
 \ov{F_j(\Gz_*,u_*)}= F_j(\Gz,u), \quad (\Gz,u)\in\CR,
 \label{2.32.3}
 \eeq
 where $(\Gz_*,u_*)=(\bar\Gz,-u(\bar\Gz))$ is the point symmetrical to the point $(\Gz,u)$ with respect to
 the line along which the two sheets $\Bbb C_1$ and $\Bbb C_2$ of the surface are connected.
 Note that if $(\Gz,u)\in \Bbb C_1$, then  $(\Gz_*,u_*)\in \Bbb C_2$.
 
 At the points $(\Gz_\infty,u_\infty)\in{\Bbb C_1}$ and  $(\Gz^*_\infty,u^*_\infty)\in{\Bbb C_2}$
 ($\Gz_\infty^*=\bar\Gz_\infty$), both of the functions $F_1(\Gz,u)$
 and $F_2(\Gz,u)$ have order-2 poles with zero residues. At the branch points of the surface, they have simple poles   
 (in the sense of Riemann surfaces). At the two infinite points of the surface
they have order-2 zeros.  
The  Schwarz boundary conditions   (\ref{2.7}) can be rewritten as the two symmetric Riemann-Hilbert problems
on the elliptic surface $\CR$ 
\beq
F_j^+(\Gx,u)-F^-_j(\Gx)=0,\quad  \Gx\in l=l_0\cup l_1.
\label{2.22}
\eeq
The solutions of these problems are rational function in the surface $\CR$. When written in the first sheet, they have the form
$$
F_-(\Gz)=\fr{1}{(\Gz-\Gz_\infty)^2}\left(
A_0^-+\fr{i(A_1^-+A_2^-\Gz+A_3^-\Gz^2)}{f(\Gz)}
\right), \quad \Gz\in{\Bbb C}\setminus l,
$$
\beq
F_+(\Gz)=\fr{1}{(\Gz-\Gz_\infty)^2}\left(
iA_0^++\fr{A_1^++A_2^+\Gz+A_3^+\Gz^2}{f(\Gz)}
\right), \quad \Gz\in{\Bbb C}\setminus l,
\label{2.21}
\eeq
where $A_j^\pm$ are real constants to be determined.
Denote 
\beq
c_{-1}=c'+ic'', \quad b\pm \bar a=\Ga^\pm+i\Gb^\pm.
\label{2.23}
\eeq
Due to (\ref{2.4}) the four parameters $\Ga^\pm$ and $\Gb^\pm$ are expressed through the 
loading data as
\beq
\Ga^+=\Gs_2^\infty-p, \quad  \Ga^-=p-\Gs_1^\infty,\quad \Gb^+=\tau^\infty-\tau, \quad \Gb^-=\tau^\infty+\tau. 
\label{2.23.0}
\eeq
On expanding the functions (\ref{2.21}) in a neighborhood of the point $\Gz=\Gz_\infty$
and satisfying the first condition in (\ref{2.7'}) we 
determine $A_0^\pm$
\beq
A_0^+=-c'\Gb^+-c''\Ga^+, \quad A_0^-=-c'\Ga^-+c''\Gb^-, 
\label{2.23.1}
\eeq
and derive four real equations for the other six coefficients
$$
A_1^\pm+\Gz_\infty A_2^\pm+\Gz_\infty^2A_3^\pm=d^\pm,
$$
\beq
A_2^\pm+2\Gz_\infty A_3^\pm=
\fr{d^\pm p_2'(\Gz_\infty)}{2p_2(\Gz_\infty)},
\label{2.23.2}
\eeq
where
\beq
d^+=\sqrt{|p_2(\Gz_\infty)|}(c'\Ga^+-c''\Gb^+), \quad d^-=\sqrt{|p_2(\Gz_\infty)|}(c'\Gb^-+c''\Ga^-).
\label{2.23.3}
\eeq 
Denoting
\beq
A_0=iA_0^+-A_0^-, \quad A_j=A_j^+-iA_j^-, \quad j=1,2,3,
\label{2.23.4}
\eeq
and 
applying formula (\ref{2.8}) we obtain the derivative of the conformal map
\beq
\Go'(\Gz)=\fr{1}{2\bar a(\Gz-\Gz_\infty)^2}\left(A_0+
\fr{A_1+A_2\Gz+A_3\Gz^2}{f(\Gz)}
\right).
\label{2.24}
\eeq
In general, the map $z=\Go(\Gz)$ given by  (\ref{2.24}) is a multi-valued  function. It is
a one-to-one map if 
\beq
\int_{l_0}\Go'(\Gz)d\Gz=0, \quad \int_{l_1}\Go'(\Gz)d\Gz=0.
\label{2.25}
\eeq
The two integrals over the loops $l_0$ and $l_1$ vanish if the coefficients
$A_j^\pm$ ($j=1,2,3$) solve the four equations
$$
A_1^\pm I_0^-+A_2^\pm I_1^-+A_3^\pm I_2^-=0,
$$
\beq
A_1^\pm I_0^+-A_2^\pm I_1^++A_3^\pm I_2^+=0,
\label{2.26}
\eeq
where
\beq
I_j^\pm=\int_1^{1/k}\fr{\Gx^jd\Gx}{(\Gx\pm\Gz_\infty)^2\sqrt{|p_2(\Gx)|}}, \quad j=0,1,2.
\label{2.27}
\eeq
The system of eight equations (\ref{2.23.2}), (\ref{2.26}) for the six unknowns $A_j^\pm$ ($j=1,2,3$)
has rank 6: the third and fourth equations in (\ref{2.26}) are identically satisfied provided
$A_j^\pm$ solve  equations (\ref{2.23.2}) and the first and second equations in  (\ref{2.26}).
Upon solving the system we express the coefficients $A_j^\pm$ through the four
problem parameters $\Ga^\pm$ and $\Gb^\pm$ and the four conformal map
parameters $c'$, $c''$,  $\Gz_\infty$, and $k$ in the form
$$
A_1^\pm=\fr{d^\pm}{\Gl_0}[\Gz_\infty(\Gl_1\Gz_\infty-2)I_1^-+(1-\Gl_1\Gz_\infty)I_2^-],
$$
\beq
A_2^\pm=\fr{d^\pm}{\Gl_0}[\Gz_\infty(2-\Gl_1\Gz_\infty)I_0^-+\Gl_1 I_2^-],\quad 
A_3^\pm=-\fr{d^\pm}{\Gl_0}[(1-\Gl_1\Gz_\infty)I_0^-+\Gl_1 I_1^-].
\label{2.28}
\eeq
Here,
\beq
\Gl_0=\Gz_\infty^2 I_0^--2\Gz_\infty I_1^-+I_2^-, \quad 
\Gl_1=\fr{p'_2(\Gz_\infty)}{2p_2(\Gz_\infty)}.
\label{2.28.0}
\eeq
The map itself, in addition to the four real  parameters $c'$, $c''$, $\Gz_\infty$, 
and $k$, has an additive constant, $B$, 
\beq
\Go(\Gz)=\fr{1}{2\bar a}
\left[
-\fr{A_0}{\Gz-\Gz_\infty}
+\int_{\Gz_0}^\Gz
\fr{(A_1+A_2\Gx+A_3\Gx^2)d\Gx}
{(\Gx-\Gz_\infty)^2\sqrt{p_2(\Gx)}}
\right]+B,
\label{2.29}
\eeq
where the path of integration $\Gz_0\Gz$ does not pass through the point $\Gz_\infty\in(-1,1)$.
When a point $\Gz$ traverses the contours $l_0$ or $l_1$, the point $z=\Go(\Gz)$ traverses 
the contours $L_0$ or $L_1$, respectively.

The function $\psi(\Gz)$ may have inadmissible poles in the exterior of  $l$.
They coincide with the zeros of the derivative $\Go'(\Gz)$ or, equivalently, with 
the zeros of the function 
\beq
\Gn(\Gz)=A_1+A_2\Gz+A_3\Gz^2+A_0 p_2^{1/2}(\Gz).
\label{2.30}
\eeq
The number of inadmissible poles of the function $\psi(\Gz)$ is determined by
\beq
Z=\fr{1}{2\pi i}\left(\int_{l_0}+\int_{l_1}+\lim_{R\to\infty}\int_{\GG_R}
\right)\fr{\Gn'(\Gz)d\Gz}{\Gn(\Gz)},
\label{2.31}
\eeq
where the positive direction is chosen such that the interior of circle $\GG_R=\{|\Gz|=R\}$
and the exterior of the cuts $l$ is on the left. Noticing  that  
\beq
\Gn'(\Gz)=A_2+2A_3\Gz+\fr{A_0\Gz(2\Gz^2-1-1/k^2)}{p_2^{1/2}(\Gz)}\sim 2(A_0+A_3)\Gz, \quad \Gz\to\infty,
\label{2.31.1}
\eeq
we can transform formula (\ref{2.31}) as
\beq
Z=2+\fr{1}{2\pi i}\left(
\int_{-1/k}^{-1}-\int_1^{1/k}\right)
\left(\fr{\Gn_0^-(\Gx)}{\Gn^-(\Gx)}-\fr{\Gn_0^+(\Gx)}{\Gn^+(\Gx)}\right)\fr{d\Gx}{\sqrt{|p_2(\Gx)|}}.
\label{2.31.2}
\eeq
where 
$$
\Gn^\pm(\Gx)=\pm i \sqrt{|p_2(\Gx)|} A_0+A_1+A_2\Gx+A_3\Gx^2,
$$
\beq
\Gn_0^\pm(\Gx)=\mp i A_0\Gx(2\Gx^2-1-1/k^{2})+\sqrt{|p_2(\Gx)|} (A_2+2A_3\Gx).
\label{2.31.3}
\eeq
On evaluating the integrals in (\ref{2.31.2}) we conclude that  
if $\Gg=|b/ a|<1$, then $Z=0$, and the function $\psi(\Gz)$ is analytic in the exterior of the cuts $l=l_0\cup l_1$
and continuous up to the boundary $l$.
When $\Gg>1$, the function $\psi(\Gz)$ has four inadmissible poles; the solution does not exist.
Finally, if $\Gg=1$, then the contours $L_0$ and $L_1$ are two straight segments, and the four poles become
removable singularities of the function $\psi(\Gz)$ lying in the loops $l_0$ and $l_1$.

\begin{figure}[t]
\centerline{
\scalebox{0.5}{\includegraphics{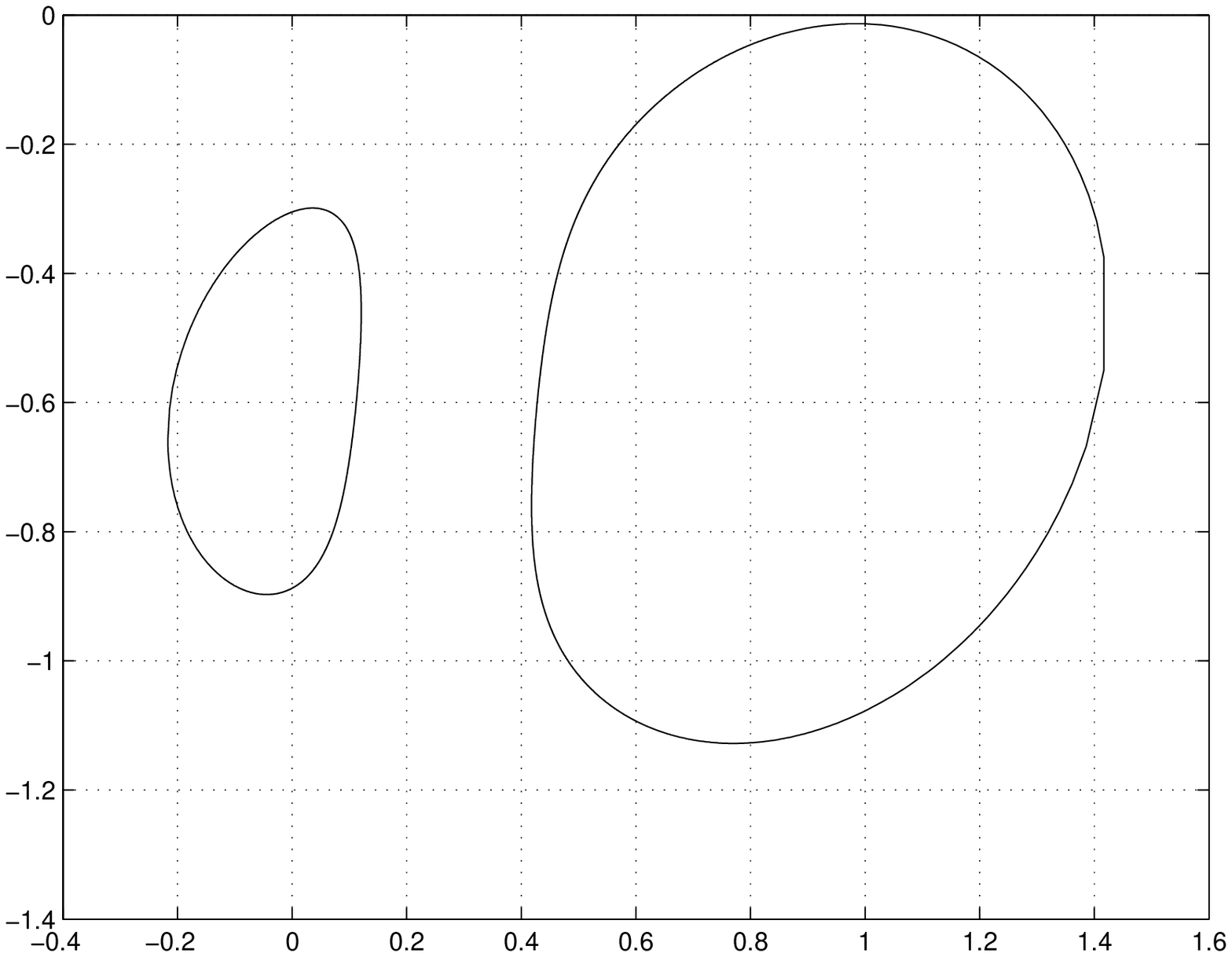}}}
\caption{
Two equal-strength holes when  $\Gz_\infty=0.5$, $p=\tau=0$,
$\Gs_1^\infty=\Gs_2^\infty=1$, $\tau^\infty=0.1$ ($\Gg=0.1$),  $k=0.01$,  $c'=1$, $c''=0$.}
\label{fig2}
\end{figure}

\begin{figure}[t]
\centerline{
\scalebox{0.5}{\includegraphics{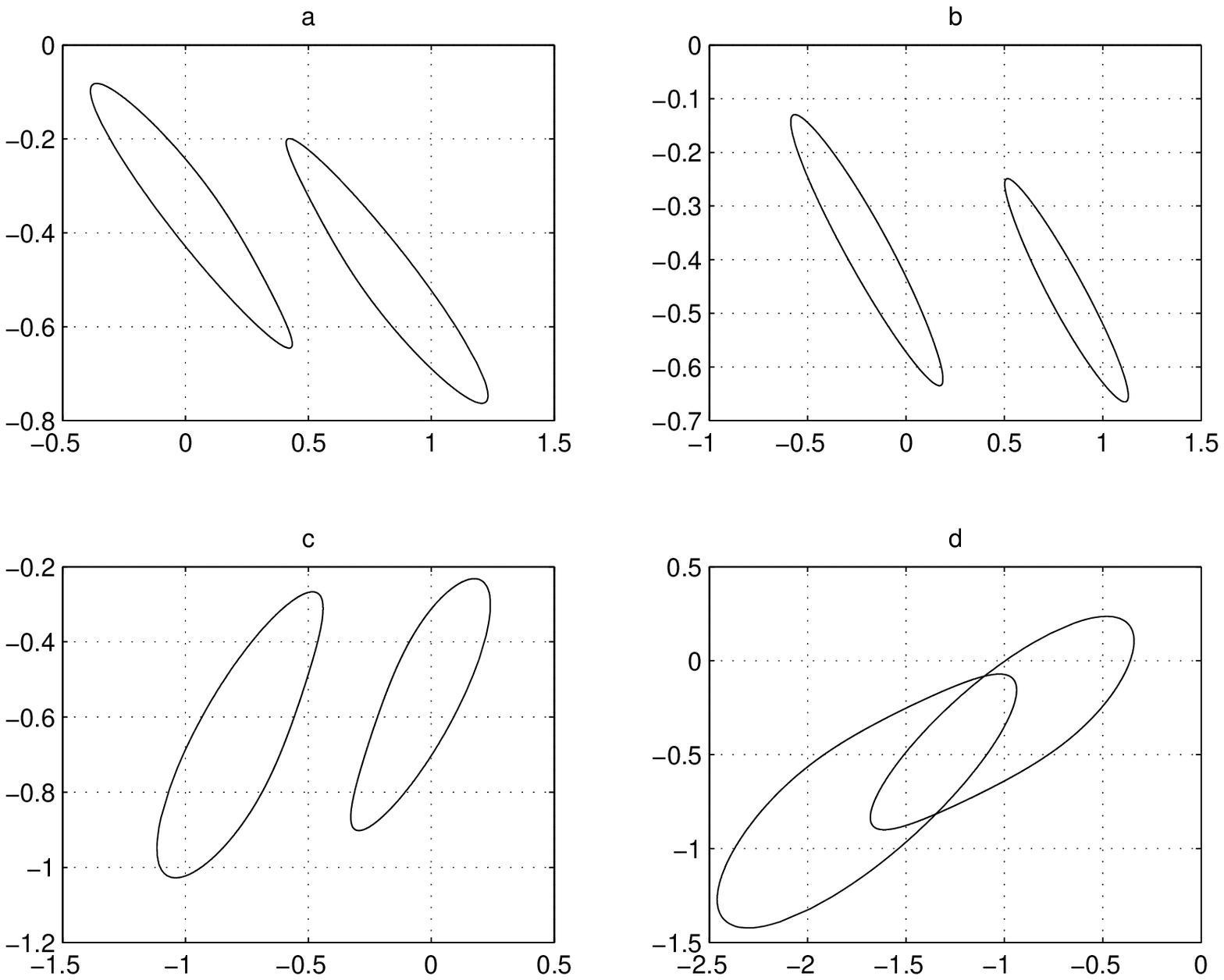}}}
\caption{Two equal-strength holes when $\Gz_\infty$ is a finite point, $p=\tau=0$, $c'=1$, and $c''=0$. 
(a): $\Gs_1^\infty=2$, $\Gs_2^\infty=1$, $\tau^\infty=-1$ ($\Gg=0.745360$),  $k=0.01$, $\Gz_\infty=0$;
(b): $\Gs_1^\infty=2$, $\Gs_2^\infty=1$, $\tau^\infty=-1$ ($\Gg=0.745360$),   $k=0.1$, $\Gz_\infty=-0.1$; 
(c): $\Gs_1^\infty=\Gs_2^\infty=2$, $\tau^\infty=1$ ($\Gg=0.5$),  $k=0.01$, $\Gz_\infty=-0.1$;
(d): $\Gs_1^\infty=\Gs_2^\infty=1$, $\tau^\infty=2$ ($\Gg=2$),  $k=0.01$, $\Gz_\infty=-0.1$.
 In case (d) $\Gg>1$,  the function
$\Psi(z)$  has four inadmissible  poles, and the solution does not exist.}
\label{fig3}
\end{figure}

\begin{figure}[t]
\centerline{
\scalebox{0.5}{\includegraphics{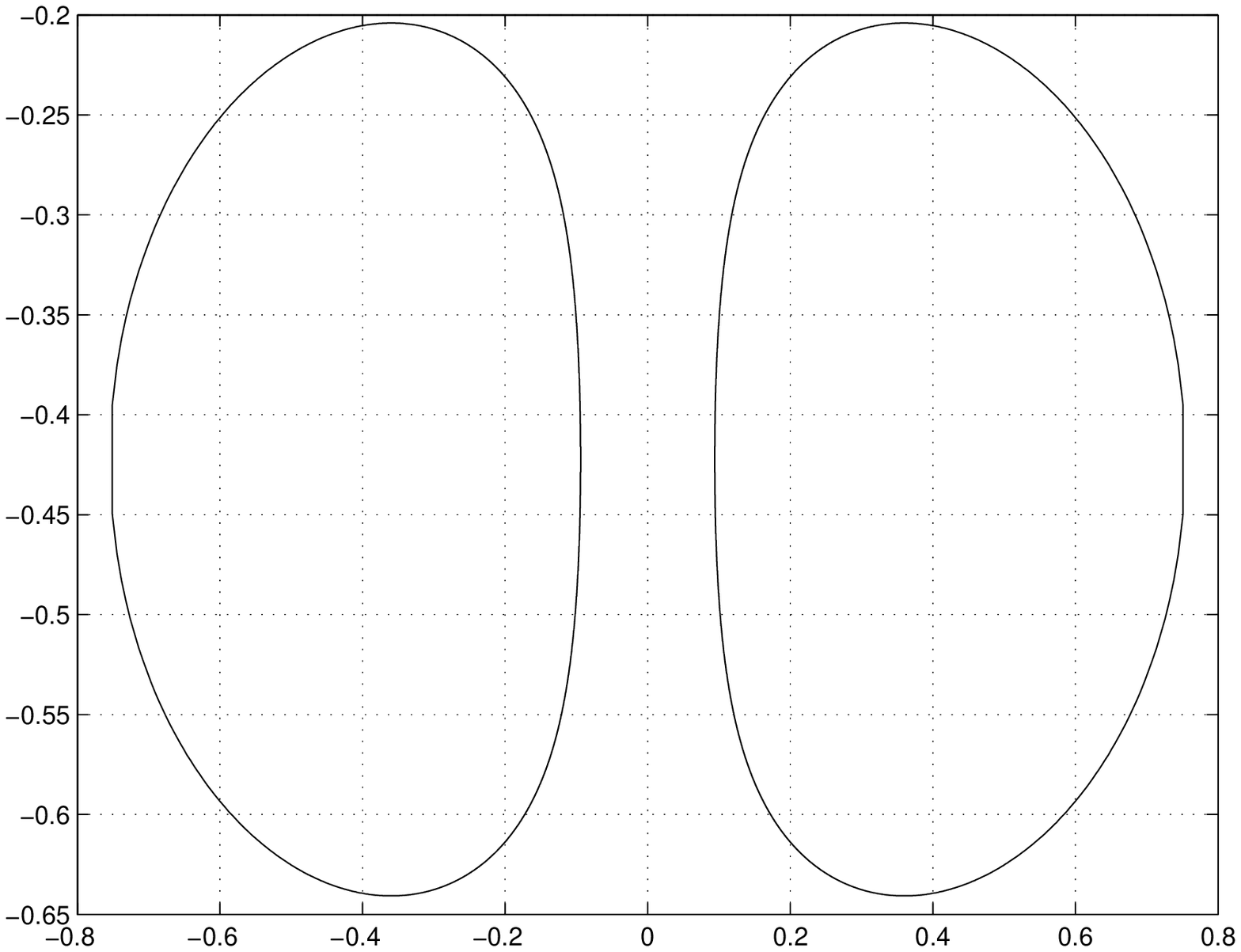}}}
\caption{
Two equal-strength symmetric holes for $\Gz_\infty=0$, $p=\tau=\tau^\infty=0$, $\Gs^\infty_1=2$, 
$\Gs^\infty_2=1$
($\Gg=1/3$), $\Gk=0.01$,
 $c'=1$, and $c''=0$.}
\label{fig4}
\end{figure} 

\begin{figure}[t]
\centerline{
\scalebox{0.5}{\includegraphics{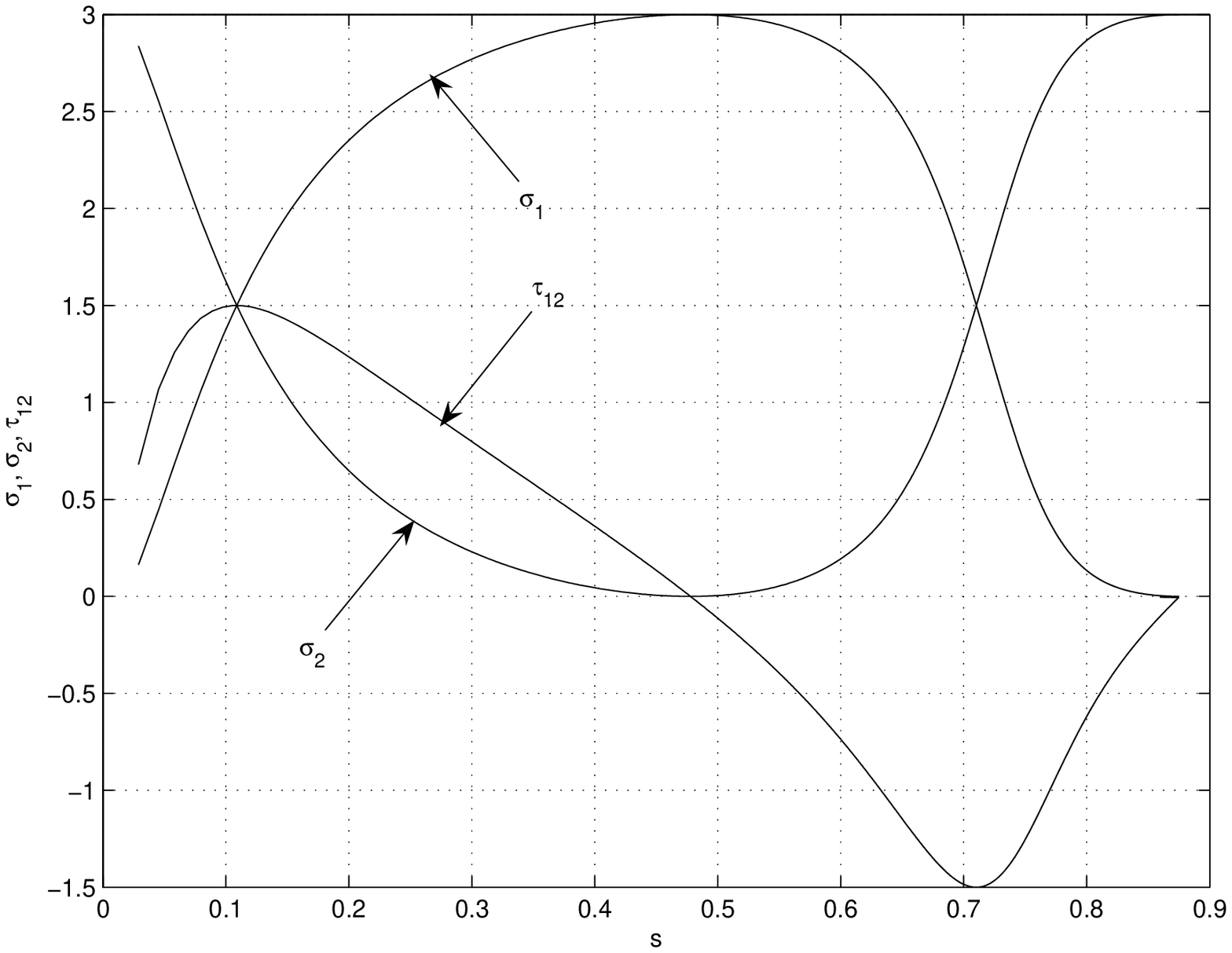}}}
\caption{The stresses $\Gs_1$, $\Gs_2$, and $\Gs_{12}$ versus the arc length $s$ on the lower half of the boundary of 
the right hole shown in Figure \ref{fig4}.
}
\label{fig5}
\end{figure}

\begin{figure}[t]
\centerline{
\scalebox{0.5}{\includegraphics{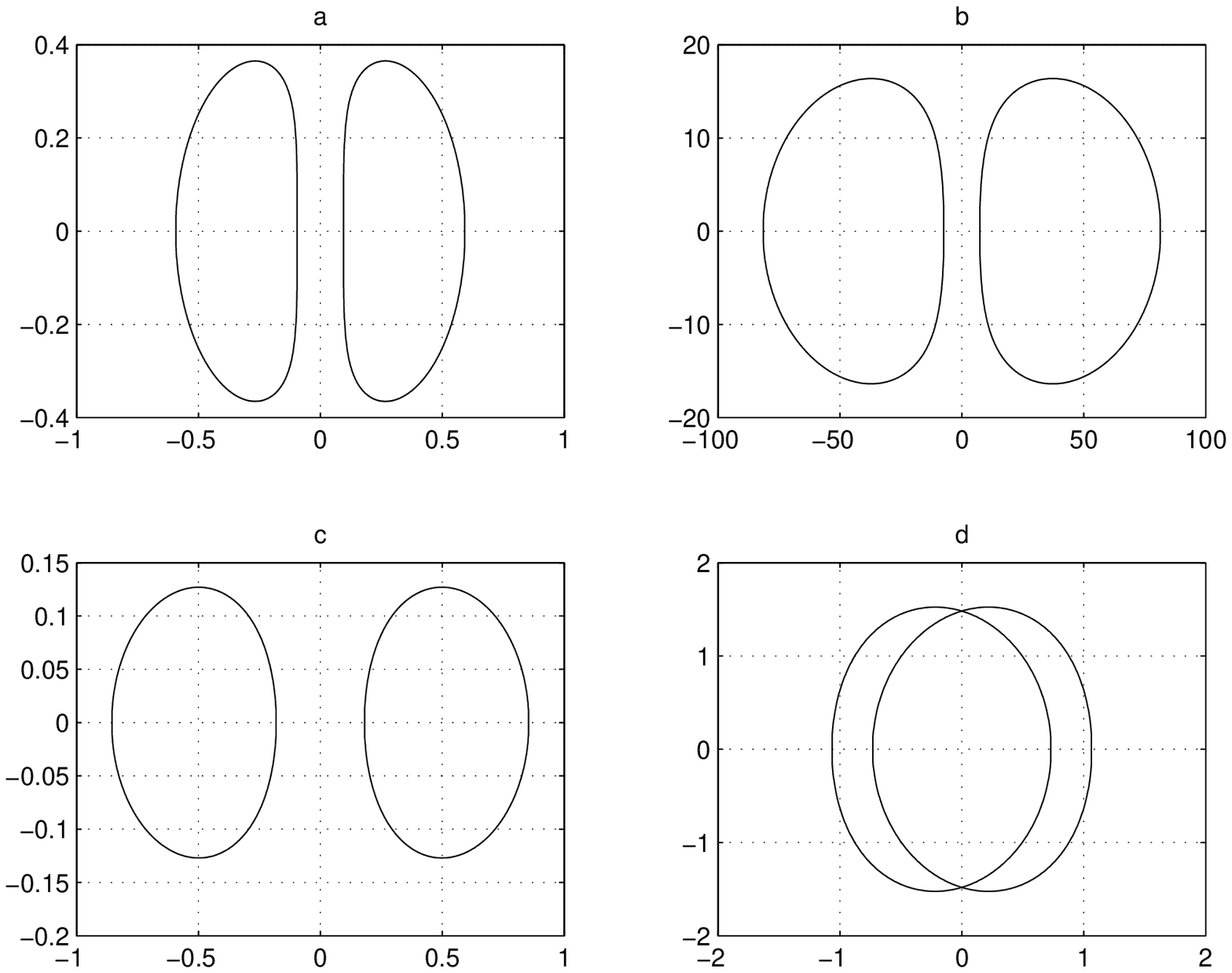}}}
\caption{
Two equal-strength symmetric holes for $\Gz_\infty=\infty$, $c'=1$, $c''=0$, $\tau^\infty=0$,
$a=1$. 
(a): $k=0.001$, $b=0$, ($\Gg=0$);
(b) $k=0.01$, $b=-0.5$ ($\Gg=0.5$);
(c) $k=0.1$, $b=-0.5$ ($\Gg=0.5$);
(d) $k=0.1$, $b=5$ ($\Gg=5$). In case (d) the function $\Psi(z)$ has four inadmissible  poles.}
\label{fig6}
\end{figure}

The profiles $L_j$ of equal-strength holes are determined from formula (\ref{2.29})
as
$$
z=\fr{1}{2\bar a}
\left[
-\fr{A_0}{\Gz-\Gz_\infty}+I_- \pm iJ(-1/k,\Gz)\right], \quad \Gz\in l_0^\pm,\quad  z\in L_0,
$$
\beq
z=\fr{1}{2\bar a}
\left[
-\fr{A_0}{\Gz-\Gz_\infty}+I_+\mp iJ(1,\Gz)\right], \quad \Gz\in l_1^\pm,\quad z\in L_1,
\label{2.31.4}
\eeq
where
\beq
I_\pm=\int_{\GG_\pm}\fr{(A_1+A_2\Gx+A_3\Gx^2)d\Gx}{(\Gx-\Gz_\infty)^2 p_2^{1/2}(\Gx)},\quad
J(d,\Gz)=\int_d^\Gz\fr{(A_1+A_2\Gx+A_3\Gx^2)d\Gx}{(\Gx-\Gz_\infty)^2 \sqrt{|p_2(\Gx)|}}.
\label{2.31.14}
\eeq
Here, $\GG_\pm$ are the segments with the staring and terminal
points $\Gz_0$ and $\pm 1$, respectively, and with no loss we assume $B=0$ and  $\Gz_0=-i$.

Some typical shapes of two nonsymmetric equal-strength holes when $\Gz_\infty$ is a finite point are shown in Figures 2 and 3. In 
Figures 2 and 3 (b) and (c), the cavities are of different area and not symmetric due to the shift of  the point $\Gz_\infty$ with respect to the center $\Gz=0$. In case 3 (a) $\Gz_\infty=0$, and the cavities have the same shape and area but not symmetric with respect to the real and imaginary axes due to the nonzero tangential stress $\tau^\infty$ applied at infinity.
In the cases shown in Figures 2 and 3 (a) -- (c) $\Gg<1$, and the function $\Psi(z)$
does not have poles in the exterior of the holes.  Referring to Figure 3(d), we observe that 
the contours $L_0$ and $L_1$ intersect each other, that is in the case $\Gg>1$
the presence of inadmissible poles 
of the function $\Psi(z)$ may not be the only one feature which indicates that the solution does
not exist. When $\Gg$ approaches $1$ and either $\Gg<1$ or $\Gg>1$, the contours become 
slim, and in the limit, when $\Gg=1$, the contours $L_0$ and $L_1$ become segments, and the function
$\Psi(z)$ is continuous everywhere in $D^e$ up to the boundary.

\subsection{Two  cavities when $\Gz_\infty$ is a finite point: the symmetric case}

The solution may be significantly simplified in the symmetric case when
$c''=0$, $\Gz_\infty=0$, $\tau=\tau^\infty=0$. Then
$$
a=\fr{\Gs-p}{2}, \quad b=\fr{\Gs_2^\infty-\Gs_1^\infty}{2},  \quad \Ga^\pm=b\pm a,\quad \Gb^\pm =0,
$$$$
d^+=\fr{c'\Ga^+}{k}, \quad d^-=0,\quad A_0^+=0, \quad
A^-_0=-c'\Ga^-,
$$
\beq
A_1^+=\fr{c'\Ga^+}{k}, \quad A_1^-=0, \quad A_2^\pm=0, \quad A^-_3=0, \quad A_3^+=-\fr{c'\Ga^+I_0^-}{kI_2^-},
\label{2.31.15}
\eeq
and the derivative of the conformal map has the form
\beq
\Go'(\Gz)=\fr{c'}{2a\Gz^2}\left[\Ga^-+\fr{\Ga^+(1-\Gz^2I_0^-/I_2^-)}{kp_2^{1/2}(\Gz)}
\right].
\label{2.31.16}
\eeq
Two sample symmetric holes are represented in Figure \ref{fig4}. For the parameters chosen, as $k\to 0$,
ellipse-like holes deform  into ``kidney"-like cavities known in the literature \cite{che}, \cite{vig1}, \cite{kan}.

It can be directly verified that $\Gs_t=\Gs$, $\Gs_n=p$, and $\tau_{nt}=\tau$, $z\in L_j$, $j=0,1$. Indeed,
on using (\ref{2.1}) to (\ref{2.4}) we have
\beq
\Gs_t=\fr{\Gs+p}{2}+\R a, \quad \Gs_n=\Gs+p-\Gs_t,\quad \tau_{tn}=\I a.
\label{2.31.16'}
\eeq
For the lower half of the boundary of the right hole shown in Figure \ref{fig4}, the variation of the stresses
$\Gs_1$, $\Gs_2$, and $\tau_{12}$ 
given by
\beq
\Gs_1=\fr12(\Gs+p)-\R\Psi(z), \quad \Gs_2=\fr12(\Gs+p)+\R\Psi(z), \quad
\tau_{12}=\I\Psi(z),
\label{2.31.16''}
\eeq
with the arc length $s$ are represented in Figure \ref{fig5}. The point $z$
traverses the contour $L_1$ in the clock wise direction, while its preimage $\Gx$ traverses the upper
side of the loop $l_1$ with the starting and terminal points $\Gx=1$ and $\Gx=1/k$, respectively. 

\subsection{Two symmetric cavities when $\Gz_\infty=\infty$}

If it is assumed that as in \cite{che}  the function $z=\Go(\Gz)$ maps the
infinite point $\Gz=\infty$ into the infinite point $z=\infty$, and the case is symmetric that is $c''=0$, $\tau=\tau^\infty=0$, then $a$ and $b$ are real,
\beq
F_-(\Gz)=c'\Ga^-, \quad F_+(\Gz)=\fr{c'\Ga^+(\Gz^2-I_2/I_0)}{p^{1/2}(\Gz)},
\label{2.31.17}
\eeq
and the function $\Go'(\Gz)$ can be represented as
\beq
\Go'(\Gz)=\fr{c'}{2a}\left[-\Ga^-+\fr{\Ga^+(\Gz^2-I_2/I_0)}{p_2^{1/2}(\Gz)}\right].
\label{2.31.18}
\eeq
 Here,
 \beq
 I_j=\int_1^{1/k}\fr{\Gx^jd\Gx}{\sqrt{|p_2(\Gx)|}}, \quad j=0,2.
 \label{2.31.19}
 \eeq
It is directly verified that the function (\ref{2.31.18}) satisfy
the conditions (\ref{2.25}), and it is a one-to-one map.
Possible poles of the function $\psi(\Gz)$  coincide with the zeros of the derivative $\Go'(\Gz)$ or, equivalently, with 
the zeros of the function 
\beq
\Gn(\Gz)=\Ga_1(\Gz^2-\Gr)-\Ga_0 p_2^{1/2}(\Gz), \quad \Ga_0=\Ga^-, \quad \Ga_1=\Ga^+, \quad \Gr=I_2/I_0.
\label{4.30}
\eeq
Due to the symmetry of the contours $l_0$ and $l_1$ the number of inadmissible poles of the function $\psi(\Gz)$ is determined by 
\beq
Z=2+\fr{1}{\pi i}\int_{1}^{1/k}\fr{\Gn^-(\Gx)[\Gn'(\Gx)]^+-\Gn^+(\Gx)[\Gn'(\Gx)]^-}
{\Gn^+(\Gx)\Gn^-(\Gx)}d\Gx,
\label{4.31}
\eeq
where $\Gn^\pm(\Gx)=\Gn(\Gx\pm i0)$, $[\Gn'(\Gx)]^\pm=\Gn'(\Gx\pm i0)$. On substituting the limiting
values  $\Gn^\pm(\Gx)$ and $[\Gn'(\Gx)]^\pm$ into the last formula we simplify it to the form
$Z=2-\Ga_0 \Ga_1 I_0/\pi$,
where we denoted
\beq
I=2\int_{1}^{1/k}
\fr{[(\Gx^2-\Gr)(2\Gx^2-1-1/k^2)+2|p_2(\Gx)|]\Gx d\Gx}
{[\Ga_1^2(\Gx^2-\Gr)^2+\Ga_0^2|p_2(\Gx)|]\sqrt{|p_2(\Gx)|}}.
\label{4.33}
\eeq
By denoting $\mu_\pm=(1/k^2\pm 1)/2$ and
making the substitutions first $\Gz^2=\mu_-\cos\Gt+\mu_+$ and then
$w=e^{i\Gt}$, we derive
\beq
I=-2i\int_{|w|=1}\fr{\{[\mu_-(w^2+1)+2\mu_+w](\mu_+-\Gr)+2w(\Gr\mu_+-1/k^2)\}dw}
{g_0(w)},
\label{2.31.8}
\eeq
where
\beq
g_0(w)=\mu_-^2(\Ga_1^2-\Ga_0^2)(w^2+1)^2+4\Ga_1^2(\mu_+-\Gr)\mu_-w(w^2+1)
+4[\Ga_1^2(\mu_+-\Gr)^2+\Ga_0^2\mu_-^2]w^2.
\label{2.31.9}
\eeq
This integral is evaluated by the theory of residues. The four zeros of the function $g_0(w)$
can be easily determined; they are
\beq
w_{1,2}=\Gd_\pm+\sqrt{\Gd_\pm^2-1}, \quad w_{3,4}=\Gd_\pm-\sqrt{\Gd_\pm^2-1}, 
\label{2.31.10}
\eeq
and
\beq
\Gd_\pm=\fr{\Ga_1^2(\Gr-\mu_+)\pm \Ga_0\sqrt{\Ga_1^2(\mu_+-\Gr)^2-\mu_-^2(\Ga_1^2-\Ga_0^2)}}
{\mu_-(\Ga_1^2-\Ga_0^2)}.
\label{2.31.11}
\eeq
The final formula for the number of inadmissible poles of the function $\psi(\Gz)$ becomes
\beq
Z=2-\Ga_*\sum_{j=1,\ldots,4; |w_j|<1}\fr{[\mu_-(w_j^2+1)+2\mu_+w_j](\mu_+-\Gr)+2w_j(\Gr\mu_+-1/k^2)}
{g_j},
\label{2.31.12}
\eeq
where $\Ga_*=\Ga_0/\Ga_1$ and
\beq
g_j=\mu_-^2(1-\Ga_*^2)w_j(w_j^2+1)+\mu_-(\mu_+-\Gr)(3w_j^2+1)+2w_j[(\mu_+-\Gr)^2+\Ga_*^2\mu_-^2].
\label{2.31.13}
\eeq
It turns out that two and only two zeros out of the four zeros $w_j$ ($j=1,2,3,4$) lie
inside the unit disc $|w|<1$.
As in the case  when $\Gz_\infty$ is a finite point in the segment (-1,1)  if $\Gg<1$,  then $Z=0$, 
 and the function $\Psi(z)$
is analytic everywhere in the domain $D^e$. If $\Gg>1$, then $Z=4$, and the function $\Psi(z)$ 
has four simple poles in the domain $D^e$.  
In the limiting case $\Gg=1$, the function $\Psi(z)$ has removable singularities in the boundary
of the domain $D^e$, and the contours $L_0$ and $L_1$ are straight segments.

Sample contours of symmetric equal-strength holes when $\Gz_\infty=\infty$ and when $\Gg<1$
are given in Figures \ref{fig6} (a) - (c). 
In Figure \ref{fig6} (d), the parameter $\Gg>1$, the function $\Psi(z)$ has four inadmissible poles, and
in addition, the contours intersect  each other; the solution does not exist.

\setcounter{equation}{0}
\section{Three cavities:  $\Gz_\infty$ is a finite point}\label{three}

Any triply connected domain $D^e$  can be considered as the image by a conformal map $z=\Go(\Gz)$
of a  parametric $\Gz$-plane  cut along three segments in the real axis, $l_0=[-1/k,-1]$,
$l_1=[k_1,k_2]$, and $l_2=[1,1/k]$, where $0<k<1$, $-1<k_1<k_2<1$. Given the domain  $D^e$ such a map
is unique. The point $z=\infty$ is the image of a certain point $\Gz_\infty=\Gz_\infty'+i\Gz_\infty''$, and, in general,
the parameters $\Gz_\infty'$, $\Gz_\infty''$, $k$, $k_1$, and $k_2$
 cannot be prescribed and should be recovered from the solution.

Let $\CR$ be the hyperelliptic  surface of the algebraic function $u^2=p_3(\Gz)$,  where
\beq
p_3(\Gz)=(\Gz^2-1)(\Gz^2-1/k^2)(\Gz-k_1)(\Gz-k_2).
\label{2.32.1}
\eeq
We fix the single branch $f(\Gz)$ of the function $p_3^{1/2}(\Gz)$ in ${\Bbb C}\setminus l$, $l=l_0\cup l_1\cup l_2$,
by the condition $f(\Gz)\sim \Gz^3$, $\Gz\to\infty$. The branch is pure imaginary on the cut sides,
$$
f(\Gz)=\pm (-1)^m i |p_3^{1/2}(\Gx)|, \quad \Gz=\Gx\pm i0, \quad \Gx\in l_m, \quad m=0,1,2,
$$
\beq
f(\Gx)=|p_3^{1/2}(\Gx)|,\quad -1<\Gx<k_1, \quad 
f(\Gx)=-|p_3^{1/2}(\Gx)|,\quad k_2<\Gx<1.
\label{2.32.2}
\eeq
Similarly to the case $n=2$, on the surface $\CR$, we introduce the functions (\ref{2.0})
which satisfy the Riemann-Hilbert problems (\ref{2.22}), (\ref{2.32.3}).  Their solution in the first sheet 
has the form
\beq
F_+(\Gz)=\fr{R_+(\Gz)}{f(\Gz)}, \quad F_-(\Gz)=\fr{iR_-(\Gz)}{f(\Gz)},
\label{2.32.4}
\eeq 
where 
$$
R_\pm(\Gz)=A_1^\pm+A_2^\pm\Gz+A_3^\pm\fr{f(\Gz)+f(\Gz_\infty)}{\Gz-\Gz_\infty}
-A_3^\pm\fr{f(\Gz)-f(\bar\Gz_\infty)}{\Gz-\bar\Gz_\infty}+(A_4^\pm+iA_5^\pm)
$$
\beq
\times\fr{f(\Gz)+f(\Gz_\infty)+f'(\Gz_\infty)(\Gz-\Gz_\infty)}{(\Gz-\Gz_\infty)^2}
-(A_4^\pm-iA_5^\pm)\fr{f(\Gz)-f(\bar\Gz_\infty)-f'(\bar\Gz_\infty)(\Gz-\bar\Gz_\infty)}{(\Gz-\bar\Gz_\infty)^2},
\label{2.32.5}
\eeq
where  $A_j^\pm$ ($j=1,2,\ldots,5$) are free real constants.
Analysis of these functions in a neighborhood of the point $\Gz_\infty$ yields
\beq
F_\pm(\Gz)=
i^{1/2\mp 1/2}\left(\fr{2(A_4^\pm+iA_5^\pm)}{(\Gz-\Gz_\infty)^2}+\fr{2A_3^\pm}{\Gz-\Gz_\infty}\right)+O(1), \quad \Gz\to\Gz_\infty.
\label{2.32.6}
\eeq
By virtue of the behavior of the functions $F_\pm(\Gz)=-c_{-1}(b\pm \bar a)(\Gz-\Gz_\infty)^{-2}+O(1)$,  
$\Gz\to\Gz_\infty$,
required, we immediately find 
$$
A_3^\pm =0,\quad A_4^+=-\fr12(c'\Ga^+-c''\Gb^+), \quad  A_4^-=-\fr12(c'\Gb^-+c''\Ga^-),
$$
\beq
A_5^+=-\fr12(c'\Gb^++c''\Ga^+), \quad  A_5^-=\fr12(c'\Ga^--c''\Gb^-).
\label{2.32.7}
\eeq
The coefficients $A_1^\pm$ and $A_2^\pm$ are still not determined in the expression
of the function $\Go'(\Gz)$ that is
$$
\Go'(\Gz)=\fr{1}{2\bar a f(\Gz)}\left[
A_1+A_2\Gz+(A_4+iA_5)\fr{f(\Gz)+f(\Gz_\infty)+f'(\Gz_\infty)(\Gz-\Gz_\infty)}{(\Gz-\Gz_\infty)^2}
\right.
$$
\beq
\left.-(A_4-iA_5)\fr{f(\Gz)-f(\bar\Gz_\infty)-f'(\bar\Gz_\infty)(\Gz-\bar\Gz_\infty)}{(\Gz-\bar\Gz_\infty)^2}\right],
\label{2.32.8}
\eeq
where $A_j=A_j^+-iA_j^-$.
To determine the coefficients $A_1$ and $A_2$, we need to guarantee that the function
$z=\Go(\Gz)$ is a single-valued map that is to force the function $\Go'(\Gz)$ to meet the three conditions
\beq
\int_{l_m} \Go'(\Gz)d\Gz=0, \quad m=0,1,2.
\label{2.32.9}
\eeq
or, equivalently, 
\beq
I_{m0}A_1+I_{m1}A_2=-J_m, \quad m=0,1,2,
\label{2.32.10}
\eeq
where
$$
I_{mj}=\int_{l_m}\fr{\Gx^j d\Gx}{f(\Gx)}, \quad j=0,1,
$$
$$
J_m=(A_4+iA_5)\int_{l_m}\fr{[f(\Gz_\infty)+f'(\Gz_\infty)(\Gx-\Gz_\infty)]d\Gx}{f(\Gx)(\Gx-\Gz_\infty)^2}
$$
\beq
+(A_4-iA_5)\int_{l_m}\fr{[f(\bar\Gz_\infty)+f'(\bar\Gz_\infty)(\Gx-\bar\Gz_\infty)]d\Gx}{f(\Gx)(\Gx-\bar\Gz_\infty)^2},
\quad m=0,1,2.
\label{2.32.11}
\eeq
The first two equations in (\ref{2.32.10}) constitute an inhomogeneous system of two complex equations with respect
to complex constants $A_1$ and $A_2$.  The coefficients of the system, the integrals $I_{mj}$
 ($j,m=0,1$) are the $A$-periods of the abelian integrals   
 \beq
 \int_{(1/k,0)}^{(\Gz,u(\Gz))}
 \fr{\Gx^j d\Gx}{u(\Gx)}, \quad j=0,1, 
 \label{2.32.12}
 \eeq
associated with the genus-2 Riemann surface $\CR$ of the algebraic function $u^2(\Gx)=p_3(\Gx)$. Therefore the $2\times 2$ matrix $\{I_{mj}\}$
($j,m=0,1$) is not singular, and the unique solution  is given by
\beq
A_1=\fr{J_1 I_{01}-J_0 I_{11}}{\GD}, \quad A_2=\fr{J_0 I_{10}-J_1 I_{00}}{\GD},
\label{2.32.13}
\eeq
where $\GD=I_{00}I_{11}-I_{01}I_{10}$.
The third equation in (\ref{2.32.10}) is transformed to the form
\beq
J_1(I_{01}I_{20}-I_{00}I_{21})+J_0(I_{10}I_{21}-I_{11}I_{20})+J_2\GD=0
\label{2.32.15}
\eeq
and satisfied identically. This is due to the fact that  the corresponding
abelian integrals in the right hand-side in (\ref{2.32.10}) can be represented as a linear combination
of the two basis integrals (\ref{2.32.12}).

Since the parameters $k$, $k_1$, $k_2$, and $\Gz_\infty=\Gz_\infty'+i\Gz_\infty''$ are free,
the derivative  $\Go'(\Gz)$ generates a five-parametric family of conformal mappings 
(we do not count the two free scaling parameters  $c_{-1}=c'+ic''$) which transforms  the slit 
domain ${\Bbb C}\setminus l$ into
the triple
connected domain $D^e$.
By integrating  the function (\ref{2.32.8}) we find the integral representation of the conformal map
$$
\Go(\Gz)=\fr{1}{2\bar a}\left\{-\fr{A_4+iA_5}{\Gz-\Gz_\infty}+\fr{A_4-iA_5}{\Gz-\bar\Gz_\infty}
+\int_{\Gz_0}^\Gz\left[
A_1+A_2\Gx+(A_4+iA_5)
\right.\right.
$$\beq
\left.\left.
\times\fr{f(\Gz_\infty)+f'(\Gz_\infty)(\Gx-\Gz_\infty)}{(\Gx-\Gz_\infty)^2}
+(A_4-iA_5)\fr{f(\bar\Gz_\infty)+f'(\bar\Gz_\infty)(\Gx-\bar\Gz_\infty)}{(\Gx-\bar\Gz_\infty)^2}\right]\fr{d\Gx}{f(\Gx)}
\right\}.
\label{2.32.16}
\eeq
To find the actual profile of the holes $L_m$, we let $\Gz$ run the cuts $l_m$ and obtain
$$
z=\CI(\Gz)\pm i\CJ(-1/k,\Gz), \quad \Gz\in l_0^\mp, \quad z\in L_0,
$$$$
z=\CI(\Gz)+\CJ(-1,k_1)\mp i\CJ(k_1,\Gz), \quad \Gz\in l_1^\mp, \quad z\in L_1,
$$
\beq
z=\CI(\Gz)+\CJ(-1,k_1)-\CJ(k_2,1)\pm i\CJ(1,\Gz), \quad \Gz\in l_2^\mp, \quad z\in L_2,
\label{2.32.17}
\eeq
where we denoted
$$
\CI(\Gz)=\fr{1}{2\bar a}\left(-\fr{A_4+iA_5}{\Gz-\Gz_\infty}+\fr{A_4-iA_5}{\Gz-\bar\Gz_\infty}\right),
$$
$$
\CJ(d,\Gz)=\fr{1}{2\bar a}
\int_{d}^\Gz\left[
A_1+A_2\Gx+(A_4+iA_5)\fr{f(\Gz_\infty)+f'(\Gz_\infty)(\Gx-\Gz_\infty)}{(\Gx-\Gz_\infty)^2}
\right.
$$\beq
\left.
+(A_4-iA_5)\fr{f(\bar\Gz_\infty)+f'(\bar\Gz_\infty)(\Gx-\bar\Gz_\infty)}{(\Gx-\bar\Gz_\infty)^2}\right]\fr{d\Gx}{|f(\Gx)|}
\label{2.32.18}
\eeq
The function $\Go'(\Gz)$ may not have zeros in the slit domain ${\Bbb C}\setminus l$. Otherwise
the functions $\psi(\Gz)$ has unacceptable  poles. 
As in the case of doubly connected domain we introduce a function $\Gn(\Gz)$ which share
the zeros with the function $\Go'(\Gz)$ and is free of singularities of $\Go'(\Gz)$,
\beq
\Go'(\Gz)=\fr{\Gn(\Gz)}{2\bar a f(\Gz) (\Gz-\Gz_\infty)^2(\Gz-\bar\Gz_\infty)^2},
\label{2.32.19}
\eeq
where
$$
\Gn(\Gz)=(A_1+A_2\Gz)(\Gz-\Gz_\infty)^2(\Gz-\bar\Gz_\infty)^2+(A_4+iA_5)(\Gz-\bar\Gz_\infty)^2[f(\Gz)+f(\Gz_\infty)
$$
\beq
+f'(\Gz_\infty)(\Gz-\Gz_\infty)]-
(A_4-iA_5)(\Gz-\Gz_\infty)^2[f(\Gz)-f(\bar\Gz_\infty)
-f'(\bar\Gz_\infty)(\Gz-\bar\Gz_\infty)].
\label{2.32.20}
\eeq
The zero counting formula applied  yields that the number of zeros, $Z$, of
the function $\Gn(\Gz)$ in the slit domain is 
\beq
Z=5+\fr{1}{2\pi i}\sum_{m=0}^2\int_{l_m}\fr{\Gn'(\Gz)d\Gz}{\Gn(\Gz)}.
\label{2.32.21}
\eeq
Here, we used the asymptotics at infinity 
$$
f(\Gz)\sim \Gz^3, \quad f'(\Gz)\sim 3\Gz^2, 
$$
\beq
\Gn(\Gz)\sim (A_2+2iA_5)\Gz^5, \quad 
\Gn'(\Gz)\sim 5(A_2+2iA_5)\Gz^4, \quad \Gz\to \infty,
\label{2.32.22}
\eeq
and the limit as $R\to\infty$ of the integral over a circle $\GG_R$ of radius $R$ centered at the origin
\beq
\lim_{R\to\infty}\int_{_R}\fr{\Gn'(\Gz)d\Gz}{\Gn(\Gz)}=5.
\label{2.32.23}
\eeq

\begin{figure}[t]
\centerline{
\scalebox{0.6}{\includegraphics{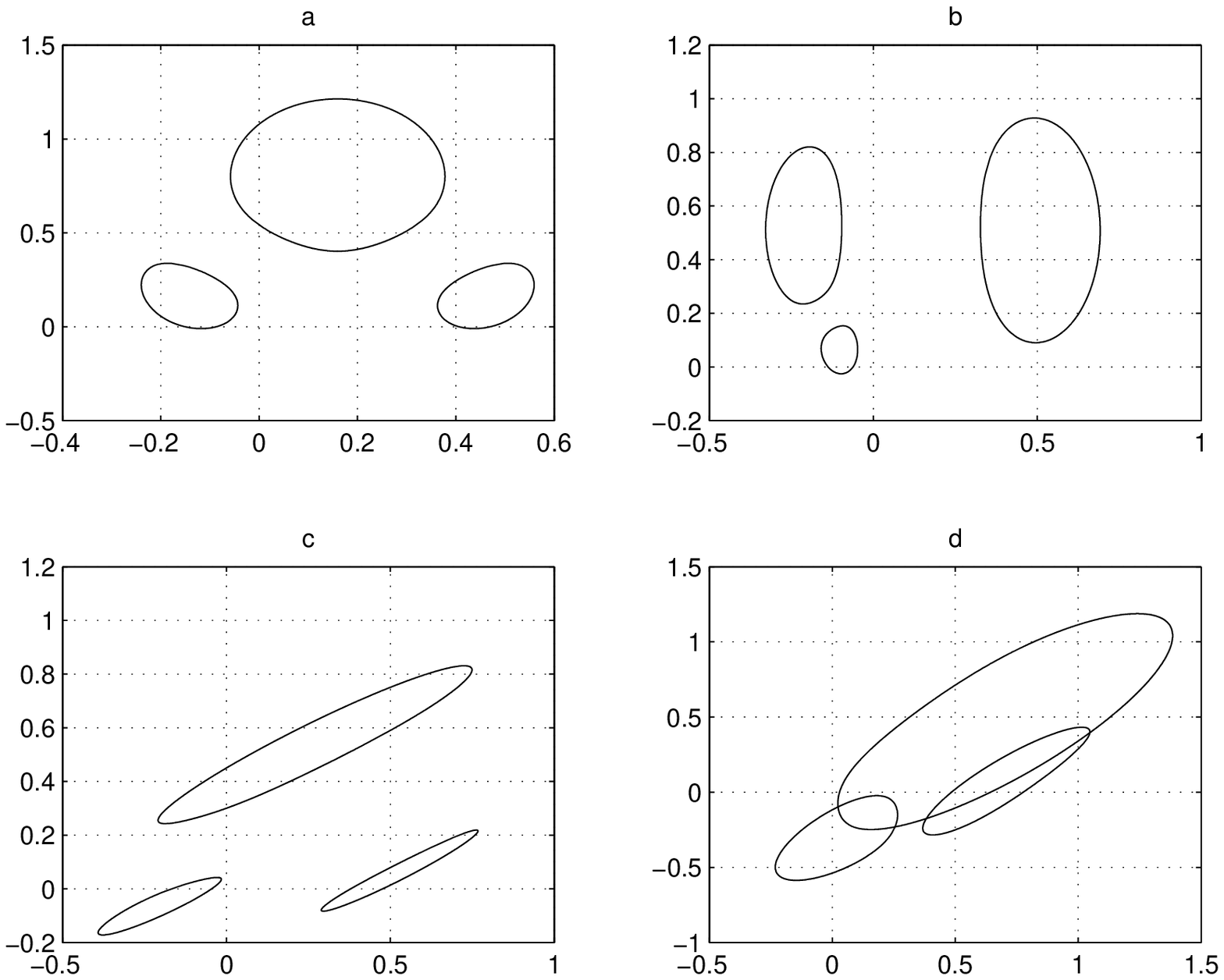}}}
\caption{
The contours $L_0$, $L_1$, and $L_2$ when $\Gz_\infty$ is a finite point, 
$c'=1, c''=0$,  $k_1=-0.8$, $k=0.2$, $k_2=0.8$, $p=\tau=0$.
(a): $\Gz_\infty=i$, $\Gs_1^\infty=1$, $\Gs_2^\infty=2$, $\tau^\infty=0$ ($\Gg=1/3$); 
(b): $\Gz_\infty=1+i$, $\Gs_1^\infty=1$, $\Gs_2^\infty=2$, $\tau^\infty=0$ ($\Gg=1/3$);
(c): $\Gz_\infty=i$, $\Gs_1^\infty=2$, $\Gs_2^\infty=1$, $\tau^\infty=1$ ($\Gg=0.7453560$); 
(d): $\Gz_\infty=i$,  $\Gs_1^\infty=\Gs_2^\infty=1$, $\tau^\infty=2$  ($\Gg=2$).
In cases (a) -- (c) the function $\Psi(z)$ has two inadmissible poles in the domain $D^e$, while in case (d) the number of 
poles is 8.}
\label{fig7}
\end{figure}

On computing the integrals in (\ref{2.32.21}) it is possible to establish that
if $\Gg=|b/a|<1$, then $Z=2$, and when $\Gg>1$, the function $\Go'(\Gz)$ has even more
zeros, $Z=8$. In the limiting case $\Gg=1$,  these points lying in the contour $l$ are removable singularities
of the function $\psi(\Gz)$
due to the relation $\ov{F_+(\Gx)+F_-(\Gx)}=-[F_+(\Gx)-F_-(\Gx)]$, $\Gx\in l$. When  $\Gg=1$, the contours
$L_j$ are straight segments in the $z$-plane.

Thus, the conformal map   which transforms any finite point $\Gz_\infty$ of the slit domain into 
the infinite point $z=\infty$  gives rise to a finite number of poles  of the function $\Psi(z)$ regardless of the value 
of the parameter $\Gg\ne 1$. 
This means that such a family of maps 
cannot be employed to identify equal-strength cavities. 
At the same time, when $\Gg<1$, the map $z=f(\Gz)$ given by (\ref{2.32.17}) 
generates some contours $L_m$ which do not intersect each other. Samples of such contours when
$\Gg<1$ and the function $\Psi(z)$ has two inadmissible poles in the domain $D^e$
are given in Figure \ref{fig7} (a) -- (c). Shown in Figure  \ref{fig7} (d) the three loops intersect each other,
$\Gg>1$, and the function $\Psi(z)$ has eight inadmissible poles in the domain $D^e$.

\setcounter{equation}{0}
\section{$n$ cavities: domain $D^e_0$, $n\ge 3$, and  $\Gz_\infty=\infty$}

The family of mappings derived in the previous section gives rise to unacceptable poles of the complex 
potential $\Psi(z)$. All the mappings share the same property, the infinite point $z=\infty$ is the image
of a certain finite point $\Gz_\infty$ in the slit domain. Since the case $\Go(\infty)=\infty$
 cannot be extracted from the solution derived in Section  \ref{three}, we consider
this case separately. Also, for generality, we assume that $n$ is not just equal to 3, but any 
finite integer $n\ge 3$. We confine ourselves to the family of domains, $D_0^e\subset D^e$, which are the images
of slit domains $\CD^e$ such that all the $n$ slits lie in the same line, $\CD^e={\Bbb C}\setminus l$,
$l=l_0\cup\ldots l_{n-1}$, and $l_j=[k_{2j},k_{2j+1}]$, $j=0,\ldots,n-1$,
$k_{2n-1}=-k_0=1$, $-1<k_1<\ldots<k_{2n-2}<1$.
The function $\Go(\Gz)$ has a simple pole  in the vicinity of the infinite point and for large $z$
can be represented
by (\ref{2.5'}). We emphasize that not every 
triply connected domain $D^e$ is the image of a slit domain $\CD^e$ such that $z=\infty$
is the image of $\Gz=\infty$. Needless to say  that not every $n$-connected  ($n\ge 4$)
domain $D^e$ is the image
of the exterior of $n$ slits lying in the same line.
On studying the family of mappings $\Go:\CD^e\to D_0^e$ when $\Go(\infty)=\infty$ 
we try to find a set  
of equal-strength holes
such that $\Go'(\Gz)$ does not have zeros in $\CD^e$, and  therefore the function $\Psi(z)$ is analytic everywhere in the domain $D_0^e$.

First we fix  the branch $f(\Gz)$ of the function $p_n^{1/2}(\Gz)$,
\beq
p_n(\Gz)=\prod_{j=0}^{2n-1} (\Gz-k_j), 
\label{6.1}
\eeq
in the domain $\CD^e$ 
by the condition $f(\Gz)\sim \Gz^n$.
The functions $F_\pm(\Gz)$ are bounded at infinity, and their counterparts defined in the Riemann 
surface have simple poles  at the branch points
of the surface $\CR$.  We have
\beq
F_+(\Gz)=\fr{1}{f(\Gz)}\sum_{j=1}^{n+1} A^+_j\Gz^{j-1}+iA^+_0, \quad 
F_-(\Gz)=\fr{i}{f(\Gz)}\sum_{j=1}^{n+1} A^-_j\Gz^{j-1}+A^-_0, 
\label{2.34}
\eeq
where     $A^\pm_j$  are arbitrary real constants. 
By expanding these functions for large $z$ and comparing these expansions with the 
asymptotics  of $F_\pm(\Gz)$ in (\ref{2.7''}) we find
$$
A^\pm_{n}=-kA^\pm_{n+1}, \quad A_{n+1}^-=\Ga^-c''+\Gb^-c',\quad A_{n+1}^+=\Ga^+c'-\Gb^+c'',
$$
\beq
A_0^-=\Ga^-c'-\Gb^-c'',\quad 
A_0^+=\Ga^+c''+\Gb^+c',
\label{2.35}
\eeq
where $k=\fr12(k_1+k_2+\ldots+k_{2n-2})$.  
The derivative of the conformal map $\Go'(\Gz)$
\beq
\Go'(\Gz)=\fr{1}{2\bar a}\left[
iA^+_0-A^-_0+\fr{1}{f(\Gz)}
\sum_{j=1}^{n+1} (A_j^+-iA_j^-)\Gz^{j-1}\right]
\label{2.36}
\eeq
has to generate a one-to-one map. This is guaranteed by the following $n$ complex conditions:
\beq
\int_{l_m} \Go'(\Gz)d\Gz=0, \quad m=0,1,\ldots,n-1.
\label{2.37}
\eeq
These conditions can be rewritten as
\beq
\sum_{j=0}^n a_{mj}(A_{j+1}^+-iA_{j+1}^-)=0, \quad m=0,1,\ldots, n-1.
\label{2.38}
\eeq
Here,
\beq
a_{mj}=\int_{l_m}\fr{\Gz^j d\Gz}{f(\Gz)}, \quad m=0,1,\ldots, n-1, \quad j=0,1,\ldots,n.
\label{2.39}
\eeq
The  integrals $a_{mj}$ ($m,j=0,1,\ldots,n-2$)  are the $A$-periods
of  the abelian integrals   
 \beq
 \int_{(1,0)}^{(\Gz,u(\Gz))}
 \fr{\Gx^j d\Gx}{u(\Gx)}, \quad j=0,1,\ldots,n-2, 
 \label{2.40}
 \eeq
associated with the genus-($n-1$) Riemann surface $\CR$ of the algebraic function $u^2(\Gx)=p_n(\Gx)$.
Therefore the matrix $a_{mj}$ ($m,j=0,1,\ldots,n-2$)  is not singular.
Denote
\beq
I_{mj}=\int_{l_m^+}\fr{\Gx^jd\Gx}{|f(\Gx)|}, \quad \quad m=0,1,\ldots, n-1, \quad  j=0,1,\ldots,n.
\label{2.41}
\eeq
The coefficients $A_j^\pm$ are uniquely determined through  the known coefficients  $A_{n+1}^\pm$ 
from the nonsingular system
\beq
\sum_{j=0}^{n-2} I_{mj} A_{j+1}^\pm=-(I_{mn}-kI_{m n-1})A_{n+1}^\pm, \quad m=0,1,\ldots,n-2.
\label{2.42.0}
\eeq
The last equation in (\ref{2.38})
is satisfied automatically because the basis of the abelian integrals (\ref{2.40}) has dimension $n-1$,
and the $n\times n$ matrix 
\beq
\left(\begin{array}{cccc}
I_{00} & \ldots & I_{0n-2} & I_{0n}-kI_{0 n-1}\\
\ldots & \ldots & \ldots  & \ldots\\
I_{n-10} & \ldots & I_{n-1n-2}  & I_{n-1n}-kI_{n-1 n-1}\\
\end{array}
\right)
\label{2.42.1}
\eeq
is singular.

To determine the number of zeros of the function $\Go'(\Gz)$, we introduce the function
\beq
\Gn(\Gz)=\sum_{j=1}^{n+1}(A^+_j-iA^-_j)\Gz^{j-1}+(iA^+_0-A^-_0)f(\Gz).
\label{2.42}
\eeq
The functions $\Gn(\Gz)$ and $\Go'(\Gz)$ share their zeros. The number of zeros of the function $\Gn(\Gz)$
coincides with the number of inadmissible poles of the function $\psi(\Gz)$ and is given by
\beq
Z=\fr{1}{2\pi i}\left(\sum_{j=0}^{n-1}\int_{l_j}+\lim_{R\to\infty}\int_{\GG_R}\right)\fr{\Gn'(\Gz)d\Gz}{\Gn(\Gz)}
=n+\fr{1}{2\pi i}\sum_{j=0}^{n-1}\int_{l_j}\fr{\Gn'(\Gz)d\Gz}{\Gn(\Gz)}.
\label{2.43}
\eeq

Our numerical tests implemented  for the case $n=3$ reveal that  if $\Gg=|b/a|<1$, then $Z=0$, and the function $\Psi(z)$
is analytic in $D_0^e$ and continuous up to the boundary.
If  $\Gg>1$, then  $Z=6$, and  the function  $\Psi(z)$ has six poles in the domain $D_0^e$. When
$|a|=|b|$, the domain $D_0^e$ is a set of straight segments, and the function   $\Psi(z)$ has removable singularities in the 
boundary of $D_0^e$. As in the simply and doubly connected cases, if $|a|<|b|$, then the problem does not have solutions.
When $|a|>|b|$, the solution exists, and  the conformal map $z=\Go(\Gz)$ is defined up to seven arbitrary constants,
an additive constant, two scaling parameters $c'$ and $c''$, and the four parameters $k_j$ ($j=1,\ldots,4$).
Integrating (\ref{2.36}) and employing (\ref{2.37}), we have
$$
z=\fr{1}{2\bar a}[(iA^+_0-A^-_0)\Gz\mp iJ(-1,\Gz)]+B,\quad \Gz\in l_0^\pm, \quad z\in L_0,
$$
$$
z=\fr{1}{2\bar a}[(iA^+_0-A^-_0)\Gz+J(k_1,k_2)\pm iJ(k_2,\Gz)]+B,\quad \Gz\in l_1^\pm, \quad z\in L_1,
$$
\beq
z=\fr{1}{2\bar a}[(iA^+_0-A^-_0)\Gz+J(k_1,k_2)-J(k_3,k_4)\mp iJ(k_4,\Gz)]+B,\quad \Gz\in l_2^\pm, \quad z\in L_2.
\label{2.45}
\eeq
Here, $B$ is an additive constant and without loss can be taken zero, and  $J$ is the real integral  
\beq
J(\Ga,\Gb)=\int_\Ga^\Gb  \fr{1}{|f(\Gx)|}\sum_{j=1}^4(A^+_j-iA^-_j)\Gx^{j-1}d\Gx.
\label{2.46}
\eeq

\begin{figure}[t]
\centerline{
\scalebox{0.6}{\includegraphics{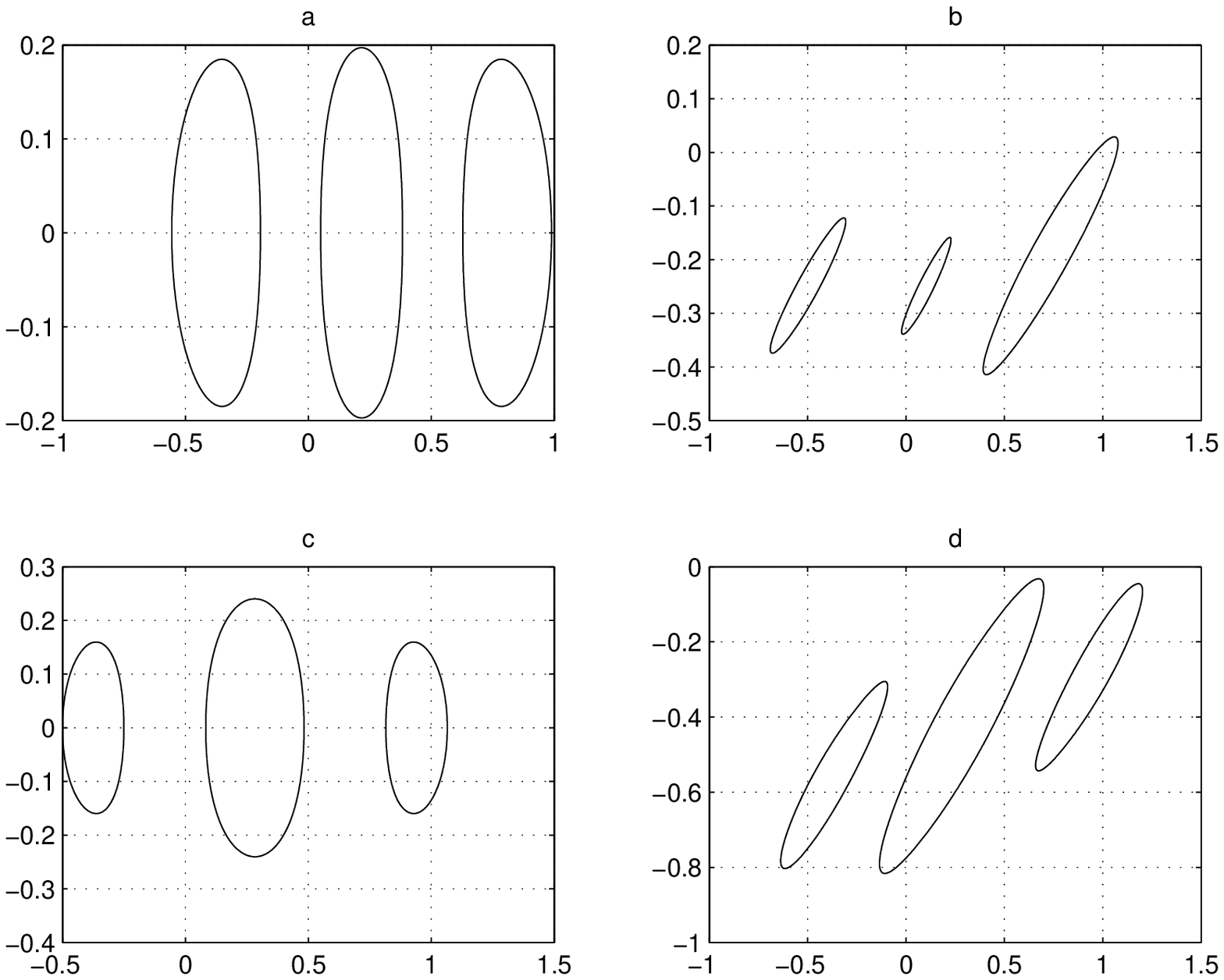}}}
\caption{
Three equal-strength holes when $\Gz_\infty=\infty$, $c'=1$, and $c''=0$, $\tau=0$.
(a): $\Gs_1^\infty=0$, $\Gs_2^\infty=1$, $\tau^\infty=0$, $p=5$ ($\Gg=1/9$),  
$k_1=-k_4=-0.35$, $k_2=-k_3=-0.3$; 
(b): $\Gs_1^\infty=2$, $\Gs_2^\infty=1$, $\tau^\infty=1$, $p=0$ ($\Gg=0.745360$),  
$k_1=-0.5$, $k_2=-0.3$, $k_3=0$, $k_4=0.1$; 
(c):  $\Gs_1^\infty=\Gs_2^\infty=1$, $\tau^\infty=0$, $p=5$ ($\Gg=0$),  
$k_1=-k_4=-0.5$, $k_2=-k_3=-0.4$; 
(d): $\Gs_1^\infty=\Gs_2^\infty=1$, $\tau^\infty=1.5$, $p=0$  ($\Gg=1.5$),  
$k_1=-k_4=-0.5$, $k_2=-k_3=-0.4$. 
In  cases (a) -- (c) the function $\Psi(z)$  is analytic, while in case (d) the function $\Psi(z)$ has
six inadmissible poles.}
\label{fig8}
\end{figure}

Figures \ref{fig8} (a)-(c)  show how the change of the loading parameter $\Gg$ and 
the conformal mapping parameters affects the profiles of equal-strength cavities in the case $n=3$
and when $\Gz_\infty=\infty$.
Figure \ref{fig8} (d) gives a sample of the contours $L_j$ when $\Gg>1$. Although the contours $L_j$ do not have common points,
the function $\Psi(z)$ has six poles in the domain $D^e$, and the solution does not exist.

\vspace{.1in}

{\bf Conclusions.} We have analyzed the inverse plane problem  of 
constructing $n$ equal-strength cavities  in an unbounded  elastic body when constant loading is
applied at infinity and to the cavities boundary. By advancing the method of conformal mappings employed 
in \cite{che} for $n=1$ and $n=2$  (two symmetric holes) to general doubly and triply connected domains
we have found by quadratures a four- and seven-parametric family of mappings
and therefore a four- and seven-parametric family of two and three equal-strength cavities, respectively.
In both cases two out of four and seven free parameters, respectively,  are scaling parameters.
For the doubly connected problem, the map $\Go$ transforms a slit domain
$\CD^e$, the exterior of two slits $[-1/k,-1]$ and $[1,1/k]$ $(0<k<1)$, into the elastic domain $D^e$, the exterior
of the holes, and $\Go(\Gz_\infty)=\infty$, $\Gz_\infty\in(-1,1)$. For the triply connected problem, we 
analyzed two cases of the preimage of the infinite point, $\Gz_\infty$  is a finite point and $\Gz_\infty=\infty$.
In the former case, $\CD^e$ is the exterior of  three slits $[-1/k,-1]$, $[k_1,k_2]$, and $[1,1/k]$,
while in the second case the slits are $[-1, k_1]$, $[k_2,k_3]$, and $[k_4,1]$. 
The conformal mappings are derived in terms of elliptic integrals for $n=2$ and hyperelliptic integrals
for $n=3$. We have also analyzed the zeros of the derivative $\Go'(\Gz)$
of the conformal mapping and shown that these zeros, if exist, generate inadmissible poles
of the solution. If  $\Gg=|b/a|<1$ ($a$ and $b$ are complex loading parameters), then the Kolosov-Muskhelishvili potential $\Psi(z)$
is free of poles  when $n=1$ and $n=2$.  
In the triply connected case, if the conformal map is chosen such that
$\Go(\infty)=\infty$, then the centers of the cavities are located in the same line, and the condition $\Gg<1$ is necessary and sufficient for the solution to exist.
If $\Gg>1$, then the function $\Psi(\Gz)$ has two, four, and six inadmissible poles in the cases $n=1$,
$n=2$ and $n=3$ ($\Go(\infty)=\infty$), respectively.
 If $\Gg=1$, then the equal-strength cavities are straight segments, and the potential $\Psi(z)$
 has removable singularities in the boundary. 
 It has also been discovered that when $\Gg>1$ and big enough,  then the contours intersect each other.
 If $n=3$ and $\Gz_\infty$ is a finite point,
 then $\Psi(\Gz)$ has two  and eight poles for the cases $\Gg<1$
and $\Gg>1$, respectively, that is there does not exist a set of three equal-strength cavities
whose centers lie not in the same line. By using the Riemann-Hilbert problem on a genus-($n-1$)
surface and the theory of  abelian integrals we have also derived an integral representation in terms
of hyperelliptic integrals for a $2n$-parametric family (two of them are scaling parameters)
of conformal 
mappings for the case $n\ge 4$ and $\Go(\infty)=\infty$, and the slits lie in the same line.
We conjecture that (i) if the centers of equal-strength cavities lie in the same line and $\Go(\infty)=\infty$, then the function $\Psi(z)$ is free of poles and the solution exists when $\Gg<1$
and the function $\Psi(z)$ has $2n$ poles in the domain $D^e$ when $\Gg>1$, and (ii)
regardless of the value of the parameter $\Gg\ne 1$, if $|\Gz_\infty|<\infty$, there does not exist
a set of $n\ge 4$ equal-strength cavities whose centers are not in the same line.

\setcounter{equation}{0}

\section*{Appendix: One cavity}

\setcounter{equation}{0}

If $n=1$, then with no loss of generality
$l_0$ and the point $\Gz_\infty$ may be selected  as   $[-1,1]$ and $\infty$, respectively.
Such a map is defined up to one real parameter, and we assume that $\I c_{-1}=0$.
Denote $p_1(\Gz)=\Gz^2-1$. We fix the branch $f(\Gz)$ of $p^{1/2}(\Gz)$ in the $\Gz$-plane cut along 
$l_0$ by the condition $p_1(\Gz)\sim \Gz$, $\Gz\to\infty$. This branch is pure imaginary on the sides of the cut, 
$f(\Gx\pm i0)=\pm i|f(\Gx)|$, $-1<\Gx<1$, and real in the real axis outside the cut.
We have
$$
F_+(\Gz)=\fr{A^+_1+ A^+_2\Gz}{f(\Gz)}+iA^+_0,\quad  F_-(\Gz)=\fr{i(A^-_1+ A^-_2\Gz)}{f(\Gz)}+A^-_0,
\eqno(A.1)
$$
where     $A^\pm_j$ $(j=0,1,2)$ are arbitrary real constants. 
Due to the asymptotics (\ref{2.7''}) of the functions $F_\pm(\Gz)$ we have 
$$
A_1^\pm=0, \quad A^+_2=c_{-1}\R(b+\bar a), \quad A^-_2=c_{-1}\I(b-\bar a),
$$
$$
A^+_0=c_{-1}\I(b+\bar a), \quad A_0^-=c_{-1}\R(b-\bar a).
\eqno(A.2)
$$
Substituting these coefficients into (A.1) and then into  (\ref{2.8}) we derive
$$
\Go'(\Gz)=\fr{c_{-1}}{2}\left[m_-+\fr{m_+\Gz}{f(\Gz)}\right], \quad m_\pm=1\pm\fr{\bar b}{\bar a},\quad 
\psi(\Gz)=\bar a\fr{(a+b)\Gz -(a-b)f(\Gz)}{(\bar a+\bar b)\Gz +(\bar a-\bar b)f(\Gz)}.
\eqno(A.3)
$$
Notice that
$$
\int_{l_0}\Go'(\Gz)d\Gz=0,
\eqno(A.4)
$$
and the map is one-to-one. The conformal map $z=\Go(\Gz)$ is defined up to an additive constant $B$,
has the form \cite{che}
$$
\Go(\Gz)=\fr{c_{-1}}{2}[m_-\Gz+m_+(\Gz^2-1)^{1/2}]+B,
\eqno(A.5)
$$
and $z=\Go(\Gz)$, $\Gz\in l_0$, is a parametric equation of a family of ellipses ($c_{-1}$ is an arbitrary nonzero real 
parameter). Denote $a_1+ib_1= \fr12c_{-1}m_-$, $a_2+ib_2= \fr12c_{-1}m_+$ and put $B=0$. Then from (A.5) we may write
$$
x=a_1\Gx\mp b_2\sqrt{1-\Gx^2}, \quad
y=b_1\Gx\pm a_2\sqrt{1-\Gx^2}.
\eqno(A.6)
$$
On excluding $\sqrt{1-\Gx^2}$ we  express  $\Gx$ through $x$ and $y$
$$
\Gx=\fr{a_2x+b_2y}{a_1a_2+b_1b_2}.
\eqno(A.7)
$$
We next square $x$ and $y$ in (A.6) and employ formula (A.7). After simple algebra we obtain
a quadratic equation in $x$ and $y$
$$
(a_2^2+b_1^2)x^2+(a_1^2+b_2^2)y^2-2(a_1b_1-a_2b_2)xy=(a_1a_2+b_1b_2)^2.
\eqno(A.8)
$$
Since its discriminant $-4\GD=-4(a_1a_2+b_1b_2)^2$ is negative, equation (A.8) 
represents an ellipse.

Now, the function $\psi(\Gz)$ has to be analytic everywhere in ${\Bbb C}\setminus l_0$.  Possible singularities
of this function coincide with the zeros   of the derivative of the map or, equivalently, with the zeros of the 
function $\Gn(\Gz)=m_-\sqrt{\Gz^2-1}+m_+\Gz$. 
The number of zeros of $\Gn(\Gz)$ inside the contour $\GG$ is given by
$$
Z=\fr{1}{2\pi i}\int_\GG\fr{\Gn'(\Gz)d\Gz}{\Gn(\Gz)}.
\eqno(A.9)
$$
Computing the integral over the contour $\GG_R$, letting $R\to\infty$, and transforming the integral over
the contour $l_0$ we obtain
$$
Z=1-\fr{m_+m_-}{\pi(m_+^2-m_-^2)}\int_{-1}^1\fr{d\Gx}{\sqrt{1-\Gx^2}(\Gx^2+\Gm^2)}, \quad
 \mu^2=\fr{m_-^2}{m_+^2-m_-^2}.
\eqno(A.10)
$$
This integral can be  computed by making the subsequent substitutions $\Gx=\cos\Gf$ and 
$e^{i\Gt}=w$
and applying the theory of residues in the $w$-plane. Eventually we derive that $Z=0$ if $|b/ a|<1$ and  $Z=2$ 
if $\Gg=|b/a|>1$. These poles can be easily determined 
$$
\Gz_{1,2}=\pm\fr{i}{2}\left(\sqrt{\fr{\bar b}{\bar a}}-\sqrt{\fr{\bar a}{\bar b}}\right).
\eqno(A.11)
$$
In the case $\Gg=1$, the ellipse becomes a segment, while the poles (A.11)  of the function $\Go'(\Gz)$,
$\Gz_{1,2}=\pm\sin\fr12(\arg a-\arg b)$, lie in the segment $l_0$. Analysis of the second formula in  (A.3) 
shows that they are removable singularities of the 
function $\psi(\Gz)$.


\vspace{.1in}

\begin{thebibliography}{99}


\bibitem{ant}
{\sc Y.A. Antipov and  V.V. Silvestrov,}  
{\it Method of Riemann surfaces in the study of supercavitating flow around two hydrofoils in a channel},  Physica D,  
235 (2007), pp. 72-81. 

\bibitem{ban}
{\sc N.V. Banichuk}, {\it  Problem of optimization of the shape of a hole  in a plate under bending},
Solid Mechanics -- Izv. AN SSSR, Mekhanika Tverdogo Tela, no.3 (1977), pp. 81-88. 

%\bibitem{che0}
%{\sc G.P. Cherepanov,} {\it Inverse problem of the theory of elasticity}, 
%Solid Mechanics -- Izv. AN SSSR, Mekhanika Tverdogo Tela, no.3 (1966), 81-82. 

\bibitem{che}
{\sc G.P. Cherepanov,} {\it Inverse problems of the plane theory of elasticity},  J. Appl. Math. Mech., 38 (1974), pp. 
915-931. 

\bibitem{cou}
{\sc R. Courant}, {\it Dirichlet's Principle, Conformal Mapping, and Minimal Surfaces}, 
Interscience Publishers, Inc., New York, N.Y., 1950.

\bibitem{esh}
{\sc J.D. Eshelby}, {\it The determination of the elastic field of an ellipsoidal inclusion, and related problems},
Proc. Roy. Soc. London A, 241 (1957), pp. 376-396.

\bibitem{gra}
{\sc Y. Grabovsky and R.V. Kohn}, {\it Microstructures minimizing the energy of a two phase elastic composite in two space dimensions. II.: The Vigdergauz microstructure}, J. Mech. Phys. Solids, 43 (1995), pp. 949- 972.

\bibitem{kan}
{\sc H. Kang, E. Kim and G. W. Milton}, {\it Inclusion pairs satisfying Eshelby's uniformity property},
SIAM, J. Appl. Math. 69 (2008), pp. 577-595.

\bibitem{kel}
{\sc M.V. Keldysh}, 
{\it Conformal mappings of multiply connected domains on canonical domains},
Uspekhi Matem. Nauk 6 (1939), pp. 90-119.

\bibitem{kos}
{\sc A.S. Kosmodamianskii}, {\it Stress state of anisotropic media with holes and cavities}, Viszcza shkola, Kiev, 1976.

\bibitem{mus}
{\sc N.I.  Muskhelishvili,} {\it Some Basic Problems of the Mathematical Theory of Elasticity}, P. Noordhoff, Ltd., Groningen, 1963.

\bibitem{neu}
{\sc H. Neuber}, {\it Zur Optimierung der Spannungskonzentration}, 
in Continuum Mechanics and Related Problems of Analysis, Nauka, Moscow (1972), pp. 375-380.

\bibitem{ru}
{\sc C.-Q. Ru and P. Schiavone}, {On the elliptic inclusion in anti-plane shear}, Math. Mech. Solids, 1
(1996), pp. 327-333.

\bibitem{sav}
{\sc G.N. Savin}, {\it Stress Distribution Around Holes}, Naukova Dumka, Kiev, 1968 (NASA Tech. Trans., Washington D.C., 1970).

\bibitem{sen}
{\sc G.P. Sendeckyj}, {\it Elastic inclusion problems in plane elastostatics}, Int. J. Solids Structures 6 (1970),
pp. 1535-1543.

\bibitem{shi}
{\sc R.-J. Shih and L.T. Wheeler}, {\it Two-dimensional inhomogeneities of minimum stress concentration},
Quart. Appl. Math. 64 (1986), pp. 567-582.

\bibitem{spr}
{\sc  G. Springer}, {\it Introduction to Riemann Surfaces}, Addison--Wesley, Reading, MA, 1956.

\bibitem{vig1} {\sc S.B. Vigdergauz}, {\it Integral equation of the inverse problem of the plane theory of elasticity},
J. Appl. Math. Mech.,  40 (1976),  pp. 518-522.

\bibitem{vig2} {\sc S.B. Vigdergauz}, {\it Effective elastic parameters of a plate with a regular system of equal-strength holes},
Solid Mechanics -- Izv. AN SSSR, Mekhanika Tverdogo Tela, 21, no.2 (1986),  pp.165-169. 

\bibitem{vig3} {\sc S.B. Vigdergauz}, {\it Constant-stress inclusions in an elastic plate},  Math. Mech. Solids, 5
(2000), pp. 265-279.

\bibitem{vig4} {\sc S.B. Vigdergauz}, {\it Stress smoothing holes in planar domains}, J. Mech. Materials Structures, 5 (2010), pp. 987-1006.

\end{thebibliography}
\end{document}